\documentclass[twocolumn]{IEEEtran}
\IEEEoverridecommandlockouts
\usepackage{color}
\usepackage{paralist}
\usepackage{booktabs, multirow, makecell}
\usepackage{subfigure, bm, comment}
\usepackage{amsmath,amsfonts,amssymb}
\usepackage{mathabx,fixmath,algorithm, algorithmic}
\usepackage{pstool}
\usepackage{hyperref}
\usepackage{cite}
\def\x{\mathbold{x}}
\def\A{\mathbold{A}}
\def\B{\mathbf{B}}

\def\Q{\mathbf{Q}}

\def\v{\mathbold{v}}
\def\c{\mathbold{c}}
\def\p{\mathbold{p}}

\urldef{\myfirstpaper}\url{http://ens.ewi.tudelft.nl/~asimonetto/TimeVarying_part1.pdf}
\urldef{\mysecondpaper}\url{http://ens.ewi.tudelft.nl/~asimonetto/TimeVarying_part2.pdf}
\def\appH{{\mathbf{H}}_{k, (K)}^{-1}}

\def\b{\mathbold{b}}

\newcommand{\uu}{\mathbold{u}}
\newcommand{\h}{\mathbold{h}}

\def\y{\mathbold{y}}

\def\H{\mathbf{H}}

\def\D{\mathbf{D}}\def\X{\mathbf{X}}

\def\transp{\mathsf{T}}

\def\qed{\hfill $\blacksquare$ }

\usepackage{theorem}
\newtheorem{remark}{Remark}

\newtheorem{proposition}{Proposition}
\newtheorem{assumption}{Assumption}

\newtheorem{theorem}{Theorem}

\newtheorem{example}{\bf Example}
\input{mysymbol.sty}
\usepackage{needspace}

\title{Decentralized Prediction-Correction Methods for Networked Time-Varying Convex Optimization}
\author{Andrea~Simonetto$^*$,~Alec~Koppel$^\dagger$,~Aryan~Mokhtari$^\dagger$,~Geert~Leus$^\S$, and~Alejandro~Ribeiro$^\dagger$

\thanks{The work in this paper is supported by NSF CAREER CCF-0952867, ONR N00014-12-1-0997, ASEE SMART, and ARL MAST CTA. This paper expands the results and presents proofs that are referenced in~\cite{Paper1Camsap},~\cite{Paper2ECC}. }
\thanks{
$^*$Andrea Simonetto is with the ICTEAM institute, Universit\'e catholique de Louvain, Belgium. Email: andrea.simonetto@uclouvain.be.}
\thanks{$^\S$ Geert Leus is with the Department of EEMCS, Delft University of Technology, The Netherlands. Email: g.j.t.leus@tudelft.nl. }
\thanks{
$^\dagger$ Alec Koppel, Aryan Mokhtari, and Alejandro Ribeiro are with the Department of ESE, University of Pennsylvania, Philadelphia, PA, USA. Emails: \{akoppel, aryanm, aribeiro\}@seas.upenn.edu. }
}
 
\begin{document}

\maketitle

\begin{abstract}
We develop algorithms that find and track the optimal solution trajectory of time-varying convex optimization problems which consist of local and network-related objectives. The algorithms are derived from the prediction-correction methodology, which corresponds to a strategy where the time-varying problem is sampled at discrete time instances and then a sequence is generated via alternatively executing predictions on how the optimizers at the next time sample are changing and corrections on how they actually have changed. Prediction is based on how the optimality conditions evolve in time, while correction is based on a gradient or Newton method, leading to Decentralized Prediction-Correction Gradient (DPC-G) and Decentralized Prediction-Correction Newton (DPC-N). We extend these methods to cases where the knowledge on how the optimization programs are changing in time is only approximate and propose Decentralized Approximate Prediction-Correction Gradient (DAPC-G) and Decentralized Approximate Prediction-Correction Newton (DAPC-N). %These methods use a first-order backward approximation to estimate the time variation of the functions. 
Convergence properties of all the proposed methods are studied and empirical performance is shown on an application of a resource allocation problem in a wireless network. We observe that the proposed methods outperform existing running algorithms by orders of magnitude. The numerical results showcase a trade-off between convergence accuracy, sampling period, and network communications.

\end{abstract}

% !TEX root = TimeVarying_part2.tex
\section{Introduction}
%mathematical description
Decentralized tracking methods are used to solve problems in which distinct agents of a network aim at minimizing a global objective that varies continuously in time. We focus on a special case of this problem, where the objective may be decomposed into two parts: the first part  is a sum of functions which are locally available at each node; the second is defined along the edges of the network, and is often defined by the cost of communication among the agents.
%applications
Problems of this kind arise, e.g., in estimation, control, and robotics \cite{Alriksson2006, Farina2010,Ogren2004, Borrelli2008,Arrichiello2006, Kim2006, Graham2012}. %The latter two settings generally amount to generating a sequence of control gains such that a network of sensors or robots may track continuously moving targets in a decentralized  manner.

One approach to continuous-time optimization problems of this kind is to sample the objective function at discrete time instances $t_k$, $k = 0,1,2,\dots,$ and then solve each time-invariant instance of the problem, via classical methods such as gradient or Newton descent. If the sampling period $h:= t_{k+1} - t_k$ is chosen arbitrarily small, then doing so would yield the solution trajectory $\y^*(t_k)$ with arbitrary accuracy. However, solving such problems for each time sample is not a viable option in most application domains, since the computation time to obtain each optimizer exceeds the rate at which the solution trajectory changes, unless $\y^*(t)$ is approximately stationary. 

Prediction-correction algorithms~\cite{%rahili2015distributed,
Paper1}, by making use of tools of {non-stationary} optimization \cite{Polyak1987,Zavala2010,Dontchev2013}, have been developed to iteratively solve convex programs which continuously vary in time. These methods operate by predicting at time $t_{k}$ the optimal solution at the discrete time instance $t_{k+1}$ via an approximation of the variation of the objective function $F$ over this time slot. Then, this prediction is revised by executing gradient or Newton descent. However, these methods are designed only for centralized settings.
We focus on time-varying convex programs in decentralized settings, where nodes can only communicate with their neighbors. As a consequence, the prediction-correction methods suggested in~\cite{Paper1} are not directly applicable. %As in the centralized case, the computation time to find the optimum for a fixed $t_k$ frequently exceeds the rate at which the optimal trajectory changes. Moreover, this issue is compounded in the decentralized setting by the additional latency required for distinct agents of a network to communicate.

One approach to solving problems of this type are decentralized running algorithms, which run at the same time-scale as the optimization problem and dynamically react to changes in the objective function. Performance guarantees for such methods yield convergence to a neighborhood of the true optimizer $\y^*(t_k)$ on the order of the sampling period $O(h)$, despite the fact that only one round of communication is allowed per discrete time step \cite{Kamgarpour2008, Farina2010, Tu2011, Bajovic2011, Zavlanos2013, Jakubiec2013, Ling2013, Simonetto2014c}. The aforementioned works mostly consider strongly convex objectives with no constraints. Notably, \cite{Jakubiec2013} and \cite{Ling2013} describe a running dual decomposition and a running alternating direction method of multipliers (ADMM) algorithm. Notice that these methods implement only correction steps and thus cannot effectively mitigate the error from the non-stationarity of the optimizer.

In this paper, we generalize the prediction-correction methodology of \cite{Paper1} to decentralized settings such that each node of a network, after communicating with its neighbors, estimates its local component of the optimal trajectory at the discrete time instance $t_{k+1}$ from information regarding the objective at time $t_{k}$, and then corrects this local prediction at time $t_{k+1}$, via additional communications within the network. To develop this generalization, in the prediction step we truncate the Taylor series of the objective function's Hessian inverse. This approximation is necessary since the computation of the objective function's Hessian inverse, which is required for the prediction step, requires global communication. In the correction step, we use decentralized approximations of gradient descent and of Newton's method to correct the predicted solution by descending towards the optimal solution of the observed objective function. In addition, we consider cases in which the prediction of how the cost function changes in time is unavailable, and must be estimated. This time-derivative approximation is particularly useful in target tracking \cite{Paper2GlobalSip} or designing learning-based control strategies \cite{Yin2014, Guan2014}.

The main contributions of the paper are the following. 
\begin{enumerate}
\item[\emph{i)}] We develop prediction-correction algorithms for a class of time-varying networked optimization problems, which can be implemented in a distributed fashion over a network of computing and communicating nodes. The correction term is either derived from a gradient method or from a (damped) Newton step. 
\item[\emph{ii)}] In order to compute the prediction (and correction for Newton) direction, we employ a novel matrix splitting technique, for which the one developed in~\cite{mokhtari2015network1,mokhtari2015network2} is a special case (only valid for adjacency matrices). The novel methodology relies on the concept of block diagonal dominance. 
\item[\emph{iii)}] We prove convergence of all the algorithms and characterize their convergence rate. For the case of the (damped) Newton correction step, we compute the (local) convergence region and argue global convergence in case of a damped step.   
\end{enumerate}

The paper is organized as follows. In Section \ref{sec:prob}, we begin by introducing the optimization problem of interest and by providing some examples for the proposed formulation. We then derive a family of algorithms which contains four distinct methods (Section \ref{sec:algorithms}). %We call these algorithms decentralized (approximate) gradient/ Newton tracking, or DeGT, DeAGT, DeNT, and DeANT, respectively. 
We analyze their convergence properties in Section \ref{sec:convg}, establishing that the sequence of iterates generated by all these algorithms converges linearly to a bounded tracking error. We observe a trade-off in the implementation between approximation accuracy and communication cost. 
In Section~\ref{sec:num}, we numerically analyze the methods on a resource allocation problem in wireless sensor networks. Lastly, in Section \ref{sec:conclusion} we conclude\footnote{
{\bf Notation.} Vectors are written as $\y\in\reals^n$ and matrices as $\A\in\reals^{n\times n}$. $\|\cdot\|$ denotes the Euclidean norm, in the case of vectors, matrices, and tensors. The gradient of the function $f(\y; t)$ with respect to $\y$ at the point $(\y,t)$ is indicated as $\nabla_{\y} f(\y; t) \in \reals^n$, while the partial derivative of the same function w.r.t. $t$ at $(\y,t)$ is $\nabla_t f(\y; t)\in \reals$. Similarly, the notation $\nabla_{\y\y} f(\y; t) \in \reals^{n\times n}$ denotes the Hessian of $f(\y;t)$ w.r.t. $\y$ at $(\y,t)$, whereas $\nabla_{t\y} f(\y; t) \in \reals^{n}$ denotes the partial derivative of the gradient of $f(\y;t)$ w.r.t. time $t$ at $(\y,t)$, i.e. the mixed first-order partial derivative vector of the objective. Consistent notation is used for higher-order derivatives.
}.

\section{Problem Formulation}\label{sec:prob}

We consider a connected undirected graph $\mathcal{G} = (V,E)$, with vertex set $V$ containing $n$ nodes and edge set $E$ containing $m$ edges. Consider $\y^i\in \reals^p$ as the decision variable of node $i$ and $t$ as a non-negative scalar that represents time. Associated with each node $i$ are time-varying strongly convex functions $f^i(\y^i ; t): \reals^{p}\times \reals_+ \to \reals $ and $g^{i,i}(\y^i; t): \reals^{p}\times \reals_+ \to \reals $. The local functions $f^i$ may be interpreted as, e.g., the merit of a particular choice of control policy \cite{Ogren2004} or statistical model \cite{Alriksson2006}.
Moreover, associated with each edge $(i,j)\in E$ is a continuously time-varying convex function $g^{i,j}(\y^i, \y^j; t): \reals^{p}\times \reals^{p}\times\reals_+ \to \reals $. These edge-wise functions represent, e.g., the cost of communicating across the network \cite{Palomar2006}.

We focus on problems where nodes aim at cooperatively minimizing the global \emph{smooth strongly convex} cost function $F:\reals^{np}\times\reals_+\to \reals$, which can be written as the sum of locally available functions $f:\reals^{np}\times\reals_+\to \reals$, and a function $g:\reals^{np}\times\reals_+\to \reals$ induced by the network structure $\mathcal{G}$. In particular, the function $f(\y;t)$ is the sum of the locally available functions $f^i(\y^i; t)$,
\begin{equation}\label{cost_functions_1}
f(\y;t):=\sum_{i\in V} f^i(\y^i; t) \; .
\end{equation}
where we have defined $\y\in \reals^{np}$ in \eqref{cost_functions_1} as the stacking of the nodes' decision variables $\y^i$, i.e., $\y=({\y^1}^\transp;\dots;{\y^n}^\transp)^\transp$. The function $g(\y;t)$ induced by the structure of the network is the sum of locally available functions $g^{i,i}(\y^i; t) $ and the functions $g^{i,j}(\y^i, \y^j; t)$ associated to the edges of the network,
\begin{equation}\label{cost_functions_2}
g(\y;t):=\sum_{i\in V} g^{i,i}(\y^i; t) + \sum_{(i,j)\in E} g^{i,j}(\y^i, \y^j; t) \; . 
\end{equation}

Our goal is to solve the time-varying convex program
\begin{equation}\label{general_cost}
\y^*(t) \!:=\! \argmin_{\y\in \reals^{np}} F(\y;t):= f(\y;t) + g(\y;t), \, \textrm{for}\,\  t\geq 0 \; ,
\end{equation}
that is the foundation of many problems in cooperative control and network utility maximization. Our goal is to enable the nodes to determine their own component of the solution $\y^*(t)$ of \eqref{general_cost} for each time $t$ in a decentralized fashion, i.e., a protocol such that each node only requires communication with neighboring nodes. 
Notice that nodes can minimize the objective function $f(\y;t)$ independently, while minimization of the function $g(\y;t)$ requires coordination and information exchange across the network. Before developing distributed protocols to solve \eqref{general_cost}, we present a couple of examples to clarify the problem setting.
%
%\vskip2mm
%
\begin{example}[Estimation of distributed processes]\label{example1}\normalfont
We consider a network of interconnected sensors monitoring a time-varying distributed process. We represent this process by a vector-valued function $\uu(\x,t)\in\reals^p$, with $\x\in\reals^3$ being the spatial coordinate, and $t$ denoting time. We assume that the process is spatially smooth so that the value of $\uu(\x,t)$ at close-by spatial coordinates is also similar. We focus on the case that a network of $n$ sensors is deployed in a spatial region $\ccalA \subset \reals^3$. The $i$-th node acquires measurements $z^i(\x^i,t)$ which are noisy linear transformations of the true process
%
%\begin{equation}
$z^i(\x^i,t) = \h^{i\,\transp} \uu(\x^i,t) + {\eta}^i(t)$, 
%\end{equation}
%
%where observations $\z^i(\x, t)$ are observations received at node $i$ which are related to the 
where $\x^i$ is the location of the sensor $i$, $\h^i$ is its regressor, and the noise ${\eta}^i(t) \sim \ccalN(0, \sigma^i)$ is Gaussian distributed independently across time with covariance $\sigma^i$. This problem setting comes up in earth sciences~\cite{Loukas2014,Roy2014} and acoustics~\cite{Martinez2013}, but it is also relevant in robotics \cite{Grocholsky2003, Leonard2010, Graham2012}. By considering the task of learning a spatially regularized least-squares estimate $\hat{\uu}\in\reals^{np}$ of the process $\uu(\x,t)$ at different locations, we obtain the time-varying networked convex program
%%%%
\begin{equation}\label{leastsquares}
\min_{\hat{\uu}^1\in\reals^p, \dots, \hat{\uu}^n\in\reals^n}\frac{1}{2}\!\sum_{i=1}^{n}\! \| \h^{i\,\transp} \hat{\uu}^i 
 - z_i(\x^i,t) \|_{\frac{1}{\sigma^i}}^2 
\!+  \frac{\beta}{2}\!\!\sum_{j\in N^i}  \!\! w^{ij} \|\hat{\uu}^i- \hat{\uu}^j \|_2^2 \; ,
\end{equation}
%%%%
where $N^i$ denotes the neighborhood of node $i$, $\hat{\uu}^i$ is the estimated value of the process $\uu(\x,t)$ at time $t$ and location $\x^i$, the constant $\beta>0$ is a regularizer that incentivizes closely located sensors to obtain similar estimates, and the nonnegative weights $w^{ij}$ may be defined according to a function of the distance between sensors. The first term in \eqref{leastsquares} defines the estimation accuracy in terms of the squared error and is identified as a sum of functions which only depend on local information, which is a special case of \eqref{cost_functions_1}. The second term in \eqref{leastsquares}
%, which incentivizes nearby sensors to select similar estimated values, 
% 
couples the decisions of node $i$ with its neighbors $j \in N^i$, and it is of the form \eqref{cost_functions_2}. Thus \eqref{leastsquares} is an instance of \eqref{general_cost}.
\end{example}
%
%\vskip2mm
%
\begin{example}[Resource allocation problems]\label{example2}\normalfont 
Consider a resource allocation problem in a wireless sensor network\cite{Xiao2006a,Palomar2006,Ghadimi2013}. Associate with sensor $i$ a time-varying utility functions $f^i:\reals^p\times \reals_{+}$  and decision variable $\y^i\in\reals^p$ representing the resources allocated to node $i$ in a network $\mathcal{G}$ of $n$ sensors. To allocate resources in this network, one must respect channel capacity and interference constraints. These constraints may be formulated in aggregate as network-flow constraints, obtaining the time-varying resource allocation problem
%
%\begin{subequations}
\begin{align}\label{res_all_1}
\min_{\y^1\in\reals^p, \dots, \y^n\in\reals^p} \ \sum_{i\in V} f^i(\y^i; t) \ \textrm{ subject to }\quad  \A\y = \b(t) \; .
\end{align} 
%\end{subequations}
%
In \eqref{res_all_1}, $\A\in\reals^{l p\times n p}$ denotes the augmented graph edge incidence matrix. The matrix $\A$ is formed by $l\times n$ square blocks of dimension $p$. If the edge $e=(j,k)$ with $j<k$ links node $j$ to node $k$ the block $(e,j)$ is $[\A]_{ej}=\mathbf{I}_p$ and the block $[\A]_{ek}=-\mathbf{I}_p$, where $\mathbf{I}_p$ denotes the identity matrix of dimension $p$. All other blocks are identically null. Moreover, the time-varying vectors $\b(t)\in\reals^{l p}$ are induced by channel capacity and rate transmission constraints.

In many situations, especially in commercial settings where the nodes are consumer devices, one seeks to solve decentralized approximations of~\eqref{res_all_1}. One way to do so is to consider the approximate augmented Lagrangian relaxation of \eqref{res_all_1}, and solve instead
\begin{align}\label{res_all_2}
\min_{\y^1\in\reals^p, \dots, \y^n\in\reals^p}& \sum_{i\in V} f^i(\y^i; t) + \frac{1}{\beta^2}\|\A\y - \b(t) \|^2 \; ,
\end{align}
which is now unconstrained~\cite{Wan2009}. Notice that the parameter $\beta>0$, which behaves similarly to a Lagrange multiplier, tunes the approximation level and penalizes the violation of the approximated constraint $\|\A\y - \b(t) \|^2$.  Observe that the first term in \eqref{res_all_2} is precisely the same as \eqref{cost_functions_1}. Moreover, block-wise decomposition of the second term yields edge-wise expressions of the form $\|(\y_i - \y_j) - \b_i(t) \|^2$, which may be identified as the functions $g^{i,j}(\y^i, \y^j; t)$ in \eqref{cost_functions_2}.

\end{example}

\section{Algorithm Development}\label{sec:algorithms}

To solve the time-varying optimization problem in \eqref{general_cost}, the first step is sampling the continuously time-varying objective function $F(\y; t)$ at time instants $t_{k}$ with $k=0,1,2,\dots$, leading to a sequence of time-invariant convex problems
\begin{equation}\label{Discrete_problem}
\y^*(t_k)\ :=\ \argmin_{\y\in\reals^{np}}\ F(\y; t_k) \qquad k\geq 0 \; .
\end{equation}
The sequence of optimal decision variables $\y^*(t_k)$ defined in \eqref{Discrete_problem} are samples of the optimal trajectory $\y^*(t)$ defined in \eqref{general_cost}. Since solving \eqref{Discrete_problem} for each time instance $t_{k}$ is impractical even for moderately sized networks, we instead devise a method to generate a sequence of approximate optimizers for~\eqref{Discrete_problem} which eventually remains close to the true optimizer $\y^*(t_k)$ in~\eqref{Discrete_problem} up to a constant error. More formally, we seek to generate a sequence $\{\y_k\}$ for which 
\begin{equation}
\limsup_{k\to \infty} \|\y_k - \y^*(t_k)\| = \textrm{const.},
\end{equation}
and whose rate, convergence, and asymptotical error constants depend on the sampling period $h$ and the number of exchanged messages per node per time instance $k$. 

To do so, we build upon prediction-correction methods, which at the current time sample $t_k$ \emph{predict} the optimal decision variable at the next time sample $t_{k+1}$, i.e., from an {arbitrary} initial variable $\y_0$, for each time $k\geq 0$, predict a new approximate optimizer as 
\begin{equation}\label{prediction}
\y_{k+1|k} = \y_k + h\, \p_{k} \; ,
\end{equation}
where index $k$ is associated with time sample $t_k$, and similarly for $k+1$ w.r.t. $t_{k+1}$, $\p_k \in \reals^{np}$ is the prediction direction, $\y_{k+1|k}$ is the predicted variable for step $k$+$1$, and $h$ is the sampling period. Then, after observing the sampled objective function at $t_{k+1}$ we correct the predicted vector $\y_{k+1|k}$ by
\begin{equation}\label{correction}
\y_{k+1} = \y_{k+1|k} + \gamma\, \c_{k+1} \; ,
\end{equation}
for a certain correction direction $\c_{k+1}\in\reals^{np}$ which defines a descent direction, with nonnegative constant step-size $\gamma >0$. 

%We generate the prediction direction to keep the suboptimality of the predicted variable $\y_{k+1|k}$ with respect to $\y^*(t_{k+1})$ as close as possible to the suboptimality of $\y_{k}$ w.r.t. $\y^*(t_k)$, using only information available at time $t_k$. 

%The predicted variable $\y_{k+1|k}$, rather than moving in a direction which is good with respect to the optimum $\y^*(t_{k})$ at time $t_k$, is chosen to anticipate the optimum $\y^*(t_{k+1})$ at the subsequent time $t_{k+1}$.
%Then, making use of information at time $t_{k+1}$, we select the correction direction as either a gradient step or an approximate Newton step. %In this section, we develop decentralized prediction-correction schemes to solve \eqref{general_cost} iteratively by a network of interconnected nodes. 
\vspace{-3mm}
\subsection{Decentralized prediction step}\label{sec:pred}

%To compute the descent direction $\p_k$ in \eqref{prediction}, we first consider the time-varying optimality  gap defined by any suboptimal trajectory $\y(t)$. In particular, the gradient of the objective function $\nabla_{\y}F(\y; t)$ with respect to $\y$ is given by 
Solving the strongly convex time-invariant problem~\eqref{Discrete_problem} accounts in finding the unique decision variable for which
\begin{equation}
\nabla_{\y}F(\y^*(t_k); t_k) = \textbf{0}. 
\end{equation}
For any other variable $\y_k \neq \y^*(t_k)$, the gradient $\nabla_{\y}F(\y_k; t_k)$ would not be null and we can use it to quantify the suboptimality of $\y$ w.r.t. $\y^*(t_k)$. 

We design the prediction direction as the one that maintains the suboptimality level when determining $\y_{k+1|k}$ (the rationale being that when arrived at optimality, we will keep it while predicting). Formally, we wish to determine $\y_{k+1|k}$ as the vector for which
\begin{equation}\label{impo}
\nabla_{\y}F(\y_{k+1|k}; t_{k+1}) = \nabla_{\y}F(\y_k; t_k).
\end{equation}
Of course, implementing~\eqref{impo} requires information at future times $t_{k+1}$ at the present $t_{k}$, an impossibility without clairvoyance. Instead, we approximate the left-hand side by adopting a Taylor expansion, obtaining, 
\begin{multline}\label{tay}
\!\!\!\!\nabla_{\y}F(\y_k; t_k) + \nabla_{\y\y}F(\y_{k}; t_{k})(\y_{k+1|k} - \y_k) + \\ h\, \nabla_{t\y}F(\y_{k}; t_{k}) = \nabla_{\y}F(\y_k; t_k),
\end{multline}
which may be reordered so that $\y_{k+1 | k}$ is on the left-hand side, yielding
%%%
\begin{equation}\label{eq.qu}
{\y}_{k+1|k} = \y_k -  h \, [\nabla_{\y\y} F(\y_k; t_k)]^{-1} \nabla_{t\y} F(\y_k; t_k) \; .
\end{equation}
%%%
The update~\eqref{eq.qu} describes the discrete-time iso-suboptimality dynamics.% for reference, its continuous-time counterpart is
%\begin{equation}\label{eq.dyn}
%\dot{\y} = - [\nabla_{\y\y} F(\y; t)]^{-1} \nabla_{t\y} F(\y; t) \; . 
%\end{equation}
%%%
\ This prediction step \eqref{eq.qu} in principle would allow us to maintain a consistent level of sub-optimality, but our focus on decentralized methods precludes its use. This is because execution of \eqref{eq.qu} requires computing the Hessian inverse $\nabla_{\y\y} F(\y_k; t_k)^{-1}$ which is not implementable by a network due to the fact that $\nabla_{\y\y} F(\y_k; t_k)^{-1}$ is a global computation. The Hessian $\nabla_{\y\y} F(\y_k; t_k)=\nabla_{\y\y} f(\y;t) + \nabla_{\y\y} g(\y;t)$ consists of two terms: The first term $\nabla_{\y\y} f(\y;t)$ is a block diagonal matrix and the second term $\nabla_{\y\y} g(\y;t)$ is a block neighbor sparse matrix that inherits the structure of the graph. Therefore, the global objective function's Hessian $\nabla_{\y\y} F(\y; t)$ %=\nabla_{\y\y} f(\y;t) + \nabla_{\y\y} g(\y;t)$ 
has the sparsity pattern of the graph and can be computed by exchanging information with neighboring nodes. Nonetheless, the Hessian inverse, required in \eqref{eq.qu}, is not neighbor sparse and its computation requires global information.

To develop a decentralized protocol to approximately execute \eqref{eq.qu}, we generalize a recently proposed technique to approximate the Hessian inverse $[\nabla_{\y\y} F(\y_k; t_k)]^{-1}$ which operates by truncating its Taylor expansion \cite{mokhtari2015network1,mokhtari2015network2}. To do so, define $\textrm{diag}[\nabla_{\y\y} g(\y_k; t_k)]$ as the block diagonal matrix which contains the diagonal blocks of the matrix $\nabla_{\y\y} g(\y_k; t_k)$, and write the Hessian $\nabla_{\y\y} F(\y_k; t_k) $ as 
\begin{equation}
\nabla_{\y\y} F(\y_k; t_k) = \D_k-\B_k \; ,
\end{equation}
where the matrices $\D_k$ and $\B_k$ are defined as
%%%
\begin{subequations}\label{decomposition}
\begin{align}\label{D_k}
\D_k &:= \nabla_{\y\y} f(\y_k; t_k) + \textrm{diag}[\nabla_{\y\y} g(\y_k; t_k)] \; , 
\\ \B_k &:= \textrm{diag}[\nabla_{\y\y} g(\y_k; t_k)] - \nabla_{\y\y} g(\y_k; t_k) \; .\label{B_k}
\end{align}
\end{subequations}
%%%
Since $F$ is strongly convex, and by Assumption~\ref{as:2} [Cf. Section~\ref{sec:convg}], the matrix $\D_k$ is a positive definite block diagonal matrix and encodes second-order local objective information. The structure of the matrix $\B_k$ is induced by that of the graph: the diagonal blocks of $\B_k$ are null and the non-diagonal block $\B_{k}^{ij}$ is nonzero and given by $-\nabla_{\y^i\y^j} g^{i,j}(\y_k^i, \y_k^j; t_k)$ iff $i$ and $j$ are neighbors. 

Given that $\D_k$ is positive definite, we can write%the objective function's Hessian $\nabla_{\y\y} F(\y_k; t_k)$ can be alternatively written as
\begin{equation}\label{hessian1}
\nabla_{\y\y} F(\y_k; t_k)=\D_k^{1/2}({\bf I} - \D_k^{-1/2} \B_k \D_k^{-1/2})\D_k^{1/2} \; .
\end{equation}
Consider now the Taylor series $
({\bf I}-\X)^{-1}=\sum_{\tau=0}^{\infty} \X^\tau$ for $\X=\D_k^{-1/2} \B_k \D_k^{-1/2} $ to write the inverse of~\eqref{hessian1} as
%%%
\begin{equation}\label{Hessian_inverse}
\hskip-.35cm[\nabla_{\y\y} F(\y_k; t_k)]^{-1}\!\! = \!\D_k^{-1/2}\!\sum_{\tau=0}^\infty \!\left(\D_k^{-1/2}\! \B_k \D_k^{-1/2}\!\right)^\tau\!\D_k^{-1/2}\!\!,\!\!
\end{equation}
%%% 
whose convergence (as well as the fact that the eigenvalues of $X$ are strictly less then one so that the Taylor series holds) will be formally proved in Appendix~\ref{app.hessian}. We approximate the Hessian inverse $[\nabla_{\y\y} F(\y_k; t_k)]^{-1}$ in \eqref{Hessian_inverse} by its $K$-th order approximate $\appH$, which is formed by truncating the series in \eqref{Hessian_inverse} to its first $K$$+$$1$ terms as
%%%
\begin{equation}\label{eq:approx}
\vspace*{-1mm} \appH = \D_k^{-1/2}\sum_{\tau=0}^K \left(\D_k^{-1/2} \B_k \D_k^{-1/2}\right)^\tau\,\D_k^{-1/2} \; .
\end{equation}
Since the matrix $\bbD_k$ is block diagonal and $\bbB_k$ is block neighbor sparse, it follows that the $K$-th order approximate inverse $\appH$ is $K$-hop block neighbor sparse, i.e. its $ij$-th block is nonzero if there is a path between nodes $i$ and $j$ with length $K$ or smaller. Substituting the approximation in \eqref{eq:approx} into \eqref{eq.qu}, the prediction step may be written as
\begin{equation}\label{eq.prediction}
{\y}_{k+1|k}= {\y}_{k} + h\,\p_{k,(K)}, 
\end{equation}
where the approximate prediction direction $\p_{k,(K)}$ is given by
%%%
\begin{equation}\label{eq.direction}
\p_{k,(K)}:=-\appH \nabla_{t \y} F(\y_k;t_k).
\end{equation}
%%%
Although the computation of the approximate prediction direction $\p_{k,(K)}$ requires information of $K$-hop neighbors, we establish that it can be computed in a decentralized manner via $K$ communication rounds among neighboring nodes.

%%%%%%%%%%%%%%%%%%%%%%%%%%%%%%%%%%%
%%%%%%%%%%%%%%%%%%%%%%%%%%%%%%%%%%%
%%%%%       P  R  O  P  O  S  I  T  I  O  N      %%%%%%%%%%%%
%%%%%%%%%%%%%%%%%%%%%%%%%%%%%%%%%%%
%%%%%%%%%%%%%%%%%%%%%%%%%%%%%%%%%%%
\begin{proposition}\label{prop.1}
Consider the prediction step~\eqref{eq.prediction} and the approximate prediction direction $\p_{k,(K)}$ in \eqref{eq.direction}. Define $\p_{k,(K)}^i$ and $\nabla_{t \y} F^i(\y_k;t_k)$ as the $i$-th sub-vector of the vectors  $\p_{k,(K)}$ and $\nabla_{t \y} F(\y_k;t_k)$ associated with node $i$. Consider $\D_k^{ij}$ and $\B_k^{ij}$ as the $ij$-th block of the matrices $\D_k$ and $\B_k$ in \eqref{D_k} - \eqref{B_k}. If node $i$ computes for $\tau \!=\! 0,\! \dots, K-1$ the recursion 
%%%
\begin{equation}\label{eq.taud}
\p_{k,(\tau+1)}^i \!=\!- (\D_k^{ii})^{-1} \Big(\! \sum_{j\in {N}^i} \B_k^{ij} \p_{k,(\tau)}^j\! +\hskip-.5mm  \nabla_{t \y} F^i(\y_k;t_k) \Big),
	\end{equation}
%%%
\vskip-3mm
\noindent with initial condition $\p_{k,(0)}^i =-(\D_k^{ii})^{-1} \nabla_{t \y} F^i(\y_k;t_k)$, the result yields the approximate prediction direction $\p_{k,(K)}^i$.% of node $i$.
\end{proposition}

%%%%%%%%%%%%%%%%%%%%%%%%%%%%%%%%%%%
%%%%%%%%%%%%%%%%%%%%%%%%%%%%%%%%%%%
%%%%%        P       R      O       O       F       %%%%%%%%%%%%%
%%%%%%%%%%%%%%%%%%%%%%%%%%%%%%%%%%%
%%%%%%%%%%%%%%%%%%%%%%%%%%%%%%%%%%%
\begin{myproof}
By direct computation. See \cite{mokhtari2015network1}, Section III for a comparable derivation.
\end{myproof}

%%%%%%%%%%%%%%%%%%%%%%%%%%%%%%%%%%%
%%%%%%%%%%%%%%%%%%%%%%%%%%%%%%%%%%%
%%%%%%%%%%%%%%%%%%%%%%%%%%%%%%%%%%%
%%%%%%%%%%%%%%%%%%%%%%%%%%%%%%%%%%%
%%%%%%%%%%%%%%%%%%%%%%%%%%%%%%%%%%%
%%%%%       M   A   I   N        M   A   T   T   E   R      %%%%%%%%%
%%%%%%%%%%%%%%%%%%%%%%%%%%%%%%%%%%%
%%%%%%%%%%%%%%%%%%%%%%%%%%%%%%%%%%%

The recursion in \eqref{eq.taud} allows for the computation the $K$-th order approximate prediction direction $\p_{k,(K)}^i$ by $K$ rounds of exchanging information with neighboring nodes. %This process can be implemented by having access to the local and neighboring information. 
The $i$-th sub-vector of the mixed partial gradient $\nabla_{t \y} F^i(\y_k;t_k) $ associated with node $i$ is given by
%%%
\begin{align}\label{new}
\nabla_{t \y} F^i(\y_k;t_k) =
	\nabla_{t\y^i} f^i(\y_k^i; t_k) + \nabla_{t\y^i}g^{i,i}(\y_k^i; t_k)
	\nonumber \\
	\quad 
	+  \sum_{j\in N^i} \!  \nabla_{t\y^i}g^{i,j}(\y_k^i,\y_k^j; t_k).
	\end{align}
%%%
Node $i$ can compute $\nabla_{t \y} F^i(\y_k;t_k) $ by having access to the decision variables of its neighbors $\y_k^j$. In addition, according to the definition of the block diagonal matrix $\bbD_k$ in \eqref{decomposition}, its $i$-th block can be written as 
%%%
\begin{align}\label{decomposition_local}
\D_k^{ii} :=& \nabla_{\y^i\y^i} f^i(\y^i_k; t_k) + \nabla_{\y^i\y^i} g^{i,i}(\y^i_k; t_k) \nonumber\\
&\quad+\sum_{j\in N^i}\nabla_{\y^i\y^i} g^{i,j}(\y^i_k, \y^j_k; t_k),
\end{align}
%%%
which is available at node $i$, after receiving $\y^j_k$. These observations imply that the initial prediction direction $\p_{k,(0)}^i \!\!=\!(\D_k^{ii})^{-1} \nabla_{t \y} F^i(\y_k;t_k)$ can be computed locally at node $i$. Further, the blocks of the neighbor sparse matrix $\B_k$ are given by
%%%
\begin{equation}\label{decomposition_local_22}
\B_k^{ij} := - \nabla_{\y^i\y^j}g^{i,j}(\y_k^i,\y_k^j; t_k) \quad \for\  j\in {N}^i,
\end{equation}
%%%
which are available at node $i$. Therefore, node $i$ can compute the recursion in \eqref{eq.taud} by having access to the $\tau$-th level approximate prediction direction $\p_{k,(\tau)}^j$ of its neighbors $j\in N^i$. After the $K$ rounds of communication with neighboring nodes, to predict the local variable ${\y}_{k+1|k}^{i} $ at step $t_k$, node $i$ executes the local update
\begin{equation}\label{eq:local_prediction}
{\y}_{k+1|k}^{i} = \y_k^i + h \,\p_{k,(K)}^i \;.
\end{equation}
%
%to predict the local variable ${\y}_{k+1|k}^{i} $ at step $t_k$. 
Thus, the prediction step in \eqref{eq.prediction} yields a decentralized protocol summarized in Algorithm~\ref{alg.ep}.

%%%%%%%%%%%%%%%%%%%%%%%%%%%%%%%%%%%
%%%%%%%%%%%%%%%%%%%%%%%%%%%%%%%%%%%
%%%%%%%       A   L   G  O  R   I   T   H   M     %%%%%%%%%%%
%%%%%%%%%%%%%%%%%%%%%%%%%%%%%%%%%%%
%%%%%%%%%%%%%%%%%%%%%%%%%%%%%%%%%%%
\begin{algorithm}[t]
\caption{Decentralized Prediction at node $i$} \label{alg.ep}\footnotesize
\textbf{Input}: The \!local \!variable \!$\y_{k}^i$, \!the \!sampling \!period \!$h$,  \!the \!approximation \!level \!$K$.

\begin{algorithmic}[1]
%\INPUT {The local variable $\y_{k}^i$. The sampling period $h=t_{k+1} - t_k$, the approximation level $K$.}
 \STATE Compute $\displaystyle{\D_k^{ii}}$ [cf. \eqref{decomposition_local}]
 \vspace{0.5mm}
 \STATE Exchange the variable $\y^i_k$ with neighbors $j\in N^i$
 \vspace{0.5mm}
 \STATE Compute the local mixed partial gradient $\nabla_{t \y} F^i(\y_k;t_k)$ [cf. \eqref{new}]
 \vspace{0.5mm}
 \STATE Compute $\displaystyle{\B_k^{ij} := - \nabla_{\y^i\y^j}g^{i,j}(\y_k^i,\y_k^j; t_k) \quad \for\  j\in {N}^i}$ [cf. \eqref{decomposition_local_22}]
\vspace{0.5mm}
 \STATE Compute $\displaystyle{\p_{k,(0)}^i = - (\D_k^{ii})^{-1} \nabla_{t \y} F^i(\y_k;t_k)}$
 \vspace{0.5mm}
%%%%%%%%%%%%%%%%%%%%%%%%%%%%%%%%%%%
%%%%   P  R   E  D   I   C   T   I   O   N         L  O   O   P   %%%%%%
%%%%%%%%%%%%%%%%%%%%%%%%%%%%%%%%%%%
  \FOR {$\displaystyle{\tau=0,1,2,\ldots, K-1}$}
  \STATE  Exchange prediction direction $\p_{k,(\tau)}^i$ with neighbors $j\in N^i$
  \STATE Compute the recursion [cf. \eqref{eq.taud}]
    
     $\displaystyle{
\p_{k,(\tau+1)}^i = -(\D_k^{ii})^{-1} \Big(\! \sum_{j\in {N}^i } \B_k^{ij} \p_{k,(\tau)}^j\! +\hskip-.5mm  \nabla_{t \y} F^i(\y_k;t_k) \Big)  }$	\vspace{-0.75mm}
	\ENDFOR
%%%%%%%%%%%%%%%%%%%%%%%%%%%%%%%%%%%
%%%%%%%%%%%%%%%%%%%%%%%%%%%%%%%%%%%	
      \STATE Predict the next trajectory
  $\displaystyle{
{\y}_{k+1|k}^{i} = \y_k^i + h \,\p_{k,(K)}^i 
}$ [cf. \eqref{eq:local_prediction}]
\end{algorithmic}
\textbf{Output}: {The predicted variable ${\y}_{k+1|k}^{i}$}.
\end{algorithm}

%%%%%%%%%%%%%%%%%%%%%%%%%%%%%%%%%%%
%%%%%%%%%%%%%%%%%%%%%%%%%%%%%%%%%%%	
%%%%%%%%%%%%%%%%%%%%%%%%%%%%%%%%%%%
%%%%%%%%%%%%%%%%%%%%%%%%%%%%%%%%%%%	

%%%%%%%%%%%%%%%%%%%%%%%%%%%%%%%%%%%
%%%%%%%%%%%%%%%%%%%%%%%%%%%%%%%%%%%
%%%%%%%       A   L   G  O  R   I   T   H   M     %%%%%%%%%%%
%%%%%%%%%%%%%%%%%%%%%%%%%%%%%%%%%%%
%%%%%%%%%%%%%%%%%%%%%%%%%%%%%%%%%%%
\begin{algorithm}[t]
\caption{Approximate Prediction at node $i$} \label{alg.ap}\footnotesize
\textbf{Input}: The \!local \!variable \!$\y_{k}^i$, \!the \!sampling \!period \!$h$,  \!the \!approximation \!level \!$K$.

\begin{algorithmic}[1]
%\INPUT {The local variable $\y_{k}^i$. The sampling period $h=t_{k+1} - t_k$, the approximation level $K$.}
 \STATE Compute $\displaystyle{\D_k^{ii}}$ [cf. \eqref{decomposition_local}]
 \vspace{0.5mm}
 \STATE Exchange the variable $\y^i_k$ with neighbors $j\in N^i$
 \vspace{0.5mm}
 \STATE Compute the approximate local mixed partial gradient $\tilde{\nabla}_{t \y} F^i_k$ [cf. \eqref{fobd}]
 \vspace{0.5mm}
 \STATE Compute $\displaystyle{\B_k^{ij} := - \nabla_{\y^i\y^j}g^{i,j}(\y_k^i,\y_k^j; t_k) \quad \for\  j\in {N}^i}$ [cf. \eqref{decomposition_local_22}]
\vspace{0.5mm}
 \STATE Compute $\displaystyle{\tilde{\p}_{k,(0)}^i = - (\D_k^{ii})^{-1} \tilde{\nabla}_{t \y} F^i_k}$
 \vspace{0.5mm}
%%%%%%%%%%%%%%%%%%%%%%%%%%%%%%%%%%%
%%%%   P  R   E  D   I   C   T   I   O   N         L  O   O   P   %%%%%%
%%%%%%%%%%%%%%%%%%%%%%%%%%%%%%%%%%%
  \FOR {$\displaystyle{\tau=0,1,2,\ldots, K-1}$}
  \STATE  Exchange prediction direction $\tilde{\p}_{k,(\tau)}^i$ with neighbors $j\in N^i$
  \STATE Compute the recursion [cf. \eqref{eq.taud}]
    
     $\displaystyle{
\tilde{\p}_{k,(\tau+1)}^i = -(\D_k^{ii})^{-1} \Big(\! \sum_{j\in {N}^i } \B_k^{ij} \tilde{\p}_{k,(\tau)}^j\! +\hskip-.5mm \tilde{\nabla}_{t \y} F^i_k \Big)  }$	\vspace{-0.75mm}
	\ENDFOR
%%%%%%%%%%%%%%%%%%%%%%%%%%%%%%%%%%%
%%%%%%%%%%%%%%%%%%%%%%%%%%%%%%%%%%%	
      \STATE Predict the next trajectory
  $\displaystyle{
{\y}_{k+1|k}^{i} = \y_k^i + h \,\tilde{\p}_{k,(K)}^i 
}$ [cf. \eqref{eq:local_prediction}]
\end{algorithmic}
\textbf{Output}: {The predicted variable ${\y}_{k+1|k}^{i}$}.
\end{algorithm}

%%%%%%%%%%%%%%%%%%%%%%%%%%%%%%%%%%%
%%%%%%%%%%%%%%%%%%%%%%%%%%%%%%%%%%%
%%%%%   S   U   B   --   S   E   C    T   I    O   N      %%%%%%%%%
%%%%%%%%%%%%%%%%%%%%%%%%%%%%%%%%%%%
%%%%%%%%%%%%%%%%%%%%%%%%%%%%%%%%%%%
\subsection{Time derivative approximation}

In practical settings, knowledge of how the function $F$ changes in time is unavailable. This issue may be mitigated by estimating the term $\nabla_{t\y} F(\y; t)$ via a first-order backward derivative: Let $\tilde{\nabla}_{t\y}{F}_k$ be an approximate version of $\nabla_{t\y}F(\y_k; t_k)$ computed as a first-order backward derivative,
%%%
\begin{equation}\label{fobd}
\tilde{\nabla}_{t\y}{F}_k =(\nabla_{\y}F(\y_k; t_k) - \nabla_{\y}F(\y_{k}; t_{k-1}))/h \; .
\end{equation}
%%%
The approximation $\tilde{\nabla}_{t\y}{F}_k$ requires only information of the previous discrete time slot. Using \eqref{fobd}, we can approximate the prediction direction as
\begin{equation}\label{eq.apprdirection}
\tilde{\p}_{k,(K)}:= - \appH \tilde{\nabla}_{t \y} F_k \; .
\end{equation}
This may be obtained in a decentralized way via $K$ rounds of communication among neighboring nodes, which may be established as a trivial extension of Proposition~\ref{prop.1}.  Algorithm~\ref{alg.ep} may be modified to instead make use of the decentralized approximate prediction step in \eqref{eq.apprdirection}, as done in Algorithm~\ref{alg.ap}. Once we obtain this local prediction of the optimizer at the next time $t_{k+1}$, using information at the current time $t_k$, the problem \eqref{general_cost} is sampled at time $t_{k+1}$. We make use of this new information in the correction step, as discussed next.
%This gives rise to the Decentralized Approximate Prediction-Correction Gradient (DAPC-G) which uses Algorithm~\ref{alg.ap} as prediction and Algorithm~\ref{alg.gc} as correction, and Decentralized Approximate Prediction-Correction Newton (DAPC-N), which uses Algorithm~\ref{alg.ap} and Algorithm~\ref{alg.nc}. In the next section, we analyze the convergence of DPC-G and DPC-N, as well as their approximated versions DAPC-G and DAPC-N.

%%%%%%%%%%%%%%%%%%%%%%%%%%%%%%%%%%%
%%%%%%%%%%%%%%%%%%%%%%%%%%%%%%%%%%%
%%%%%   S   U   B   --   S   E   C    T   I    O   N      %%%%%%%%%
%%%%%%%%%%%%%%%%%%%%%%%%%%%%%%%%%%%
%%%%%%%%%%%%%%%%%%%%%%%%%%%%%%%%%%%
\vspace{-4mm}
\subsection{Decentralized correction step}
The predicted variable $\y_{k+1|k}$ [cf. \eqref{eq.direction}] is then corrected via~\eqref{correction} by making use of the objective at time $t_{k+1}$. 
%
%, i.e.
%%%
%\begin{equation}\label{eq:correction1}
%\c_{k+1} = - \A_{k+1}^{-1} \nabla_{\y}F(\y_{k+1|k}; t_{k+1}) \; , 
%\end{equation}
%%%
%where $\A_{k+1} \in \reals^{np\times np}$ is a nonsingular ``preconditioning'' correction matrix, while $\nabla_{\y}F(\y_{k+1|k}; t_{k+1})$ is the gradient of the cost function computed at time $t_{k+1}$ evaluated at the predicted vector $\y_{k+1|k}$. Different choices for the correction matrix $\A_{k+1}$  
%
Different correction strategies give rise to different correction updates,  whose relative merits depend on the application domain at hand. % -- see Section \ref{sec:num}. 
We present two distinct correction steps next.

%%%%%%%%%%%%%%%%%%%%%%%%%%%%%%%%%%%
%%%%%%%%%%%%%%%%%%%%%%%%%%%%%%%%%%%
%%%%%   S   U   B   --  S  U   B   --   S   E   C    T   I    O   N      %%%
%%%%%%%%%%%%%%%%%%%%%%%%%%%%%%%%%%%
%%%%%%%%%%%%%%%%%%%%%%%%%%%%%%%%%%%
\vspace{2mm}
{\it \textbf{Gradient correction step:}}
After the objective at time $t_{k+1}$ is observed, we may execute the correction step \eqref{correction} with $\c_{k+1} = - \nabla_{\y}F(\y_{k+1|k}; t_{k+1}) $, resulting in
%with the matrix $\A_{k+1} = {\bf I}$, resulting in
%%%
\begin{equation}\label{correctionstep}
\y_{k+1} = \y_{k+1|k} - \gamma \nabla_{\y}F(\y_{k+1|k}; t_{k+1})\;,
\end{equation}
%%%
which is a gradient correction step. This step is computable in a decentralized fashion since the local component of the gradient $F(\y_{k+1|k}; t_{k+1})$ at node $i$ is given by
%%%
\begin{align}\label{new2}
&\nabla_{\y} F^i(\y_{k+1|k}; t_{k+1}) =
	\nabla_{\y^i} f^i(\y_{k+1|k}^i; t_{k+1}) 	 \\
&	
\!	+\! \nabla_{\y^i}g^{i,i}(\y_{k+1|k}^i; t_{k+1})
\!+\!\!\!\!  \sum_{j\in N^i} \!  \nabla_{\y^i}g^{i,j}(\y_{k+1|k}^i,\y_{k+1|k}^j; t_{k+1}).\nonumber
	\end{align}
%%%
To implement the expression in~\eqref{new2}, node $i$ only requires access to the decision variables $\y^j_{k+1|k}$ of its neighbors $j \in N^i$. Thus, if nodes exchange their predicted variable $\y^i_{k+1|k}$ with their neighbors they can compute the local correction direction $\c^i_{k+1}$ as in \eqref{new2} and update their predicted variable $\y_{k+1|k}^i$ as
%%%
\begin{equation}\label{local_cor}
\y_{k+1}^i = \y_{k+1|k}^i + \gamma \c^i_{k+1}.
\end{equation}
%%%

%%%%%%%%%%%%%%%%%%%%%%%%%%%%%%%%%%%
%%%%%%%%%%%%%%%%%%%%%%%%%%%%%%%%%%%
%%%%%%%%%%%%%%%%%%%%%%%%%%%%%%%%%%%
\begin{algorithm}[t]
\caption{Decentralized Gradient Correction at node $i$} \label{alg.gc}\footnotesize
\textbf{Input}: The local predicted variable $\y_{k+1|k}^i$.  The step-size $\gamma$.

\begin{algorithmic}[1]
\STATE Exchange the predicted variable ${\y}_{k+1|k}^{i}$ with neighbors $j\in N^i$
\STATE Observe $F^i(\cdot ; t_{k+1})$, find $\c^i_{k+1}=-\nabla_{\y} F^i(\y_{k+1|k}; t_{k+1}) $ [cf. \eqref{new2}]
\STATE Correct the trajectory    
  $\displaystyle{
\y_{k+1}^i = \y_{k+1|k}^i + \gamma \c^i_{k+1}
}$   [cf. \eqref{local_cor}]
\end{algorithmic}
\textbf{Output}: {The corrected variable ${\y}_{k+1}^{i}$}.
\end{algorithm}
%%%%%%%%%%%%%%%%%%%%%%%%%%%%%%%%%%%
%%%%%%%%%%%%%%%%%%%%%%%%%%%%%%%%%%%
%%%%%%%%%%%%%%%%%%%%%%%%%%%%%%%%%%%

We call DPC-G as the Decentralized Prediction-Correction method that uses gradient descent in the correction step (Algorithm~\ref{alg.gc}) and the {\it exact} prediction step (Algorithm~\ref{alg.ep}) in the prediction step. We call DAPC-G as the Decentralized Approximate Prediction-Correction method that uses gradient descent in the correction step (Algorithm~\ref{alg.gc}) and the {\it approximate} prediction step (Algorithm~\ref{alg.ap}) in the prediction step. Both DPC-G and DAPC-G require $K+2$ communication rounds among neighboring nodes per time step.

\vspace{2mm}
{\it \textbf{Newton correction step:}}
%\textbf{(2)} 
The correction step in \eqref{correction} could also be considered as a Newton step if we used $\c_{k+1} = - \nabla_{\y\y}F(\y_{k+1|k}; t_{k+1})^{-1}\nabla_{\y}F(\y_{k+1|k}; t_{k+1})$. % and step-size $\gamma$ set to $1$.  
%
%
% the Hessian $ \nabla_{\y\y}F(\y_{k+1|k}; t_{k+1})$ as the preconditioning matrix, i.e., $\A_{k+1} =  \nabla_{\y\y}F(\y_{k+1|k}; t_{k+1})$. Note that the step-size $\gamma$ in~\eqref{correction} will be set to $\gamma = 1$. 
%
However, as in the discussion regarding the prediction step, computation of the Hessian inverse $ \nabla_{\y\y}F(\y_{k+1|k}; t_{k+1})^{-1}$ requires global communication. Consequently, we approximate the Hessian inverse $\nabla_{\y\y}F(\y_{k+1|k}; t_{k+1})^{-1}$ by truncating its Taylor series as in \eqref{eq:approx}. To be more precise, we define 
${\bf H}_{k+1|k,(K')}^{-1}$ as the $K'$-th level approximation of the Hessian inverse as %$\nabla_{\y\y}F(\y_{k+1|k}, t_{k+1})^{-1}$ as
%%%
\begin{equation}\label{eq:approx22}\vspace{-4mm}
{\bf H}_{k+1|k,(K')}^{-1} \!=\! \D_{k+1|k}^{-1/2}\!\!\sum_{\tau=0}^{K'} \left(\D_{k+1|k}^{-1/2} \B_{k+1|k} \D_{k+1|k}^{-1/2}\right)^\tau \!\!\D_{k+1|k}^{-1/2},
\end{equation}
where the matrices $\D_{k+1|k}$ and $\B_{k+1|k}$ are defined as 
\begin{subequations}\label{decomposition2}
\begin{align}\label{pofak}
\D_{k+1|k} &\!:=\! \nabla_{\y\y} f(\y_{k+1|k}; t_{k+1}) \!+\! \textrm{diag}[\nabla_{\y\y} g(\y_{k+1|k}; t_{k+1})] \; , 
\\ \B_{k+1|k} &\!:=\! \textrm{diag}[\nabla_{\y\y} g(\y_{k+1|k}; t_{k+1})]\! -\! \nabla_{\y\y} g(\y_{k+1|k}; t_{k+1}) \; .
\end{align}
\end{subequations}
Notice that the only difference between the decomposition matrices $\D_{k+1|k}$ and $\B_{k+1|k}$ for the correction step and the matrices $\D_{k}$ and $\B_{k}$ for the prediction step is the arguments for the inputs $\y$ and $t$. The prediction matrices $\D_{k}$ and $\B_{k}$ are evaluated for the function $F(.;t_k)$ and the variable $\y_k$, while the correction matrices are evaluated for the function $F(.;t_{k+1})$ and the variable $\y_{k+1|k}$.

Thus, we can approximate the exact Hessian inverse $\nabla_{\y\y}F(\y_{k+1|k}; t_{k+1})^{-1}$ with ${\bf H}_{k+1|k,(K')}^{-1} $ as in \eqref{eq:approx22} and apply the correction step as
%%%
\begin{equation}\label{correctionstepnewton}
\y_{k+1} = \y_{k+1|k} - \gamma\, {\bf H}_{k+1|k,(K')}^{-1}\nabla_{\y}F(\y_{k+1|k}; t_{k+1}) \; ,
\end{equation}
%%%
which requires $K'$ exchanges of information among neighboring nodes. In practice, one can use the same algorithm for the prediction direction $\p_{k,(K)}$ to compute the correction direction  $\c_{k,(K')}:=-{\bf H}_{k+1|k,(K')}^{-1}\nabla_{\y}F(\y_{k+1|k}; t_{k+1})$, where now the gradient takes the place of the time derivative.

%%%%%%%%%%%%%%%%%%%%%%%%%%%%%%%%%%%
%%%%%%%%%%%%%%%%%%%%%%%%%%%%%%%%%%%
%%%%%%%       A   L   G  O  R   I   T   H   M     %%%%%%%%%%%
%%%%%%%%%%%%%%%%%%%%%%%%%%%%%%%%%%%
%%%%%%%%%%%%%%%%%%%%%%%%%%%%%%%%%%%
\begin{algorithm}[t]
\caption{Decentralized Newton Correction at node $i$} \label{alg.nc}\footnotesize
\textbf{Input}: The local predicted variable $\y_{k+1|k}^i$.  The approximation level $K'$. The step-size $\gamma$.

\begin{algorithmic}[1]
\STATE Exchange the predicted variable ${\y}_{k+1|k}^{i}$ with neighbors $j\in N^i$
\STATE Observe $F^i(\cdot ; t_{k+1})$, compute $\nabla_{\y} F^i(\y_{k+1|k}; t_{k+1}) $ [cf. \eqref{new2}]
 \STATE Compute matrices $\displaystyle{\D_{k+1|k}^{ii}}$ and $\displaystyle{\B_{k+1|k}^{ij}}$, $\displaystyle{ j \in {N}^i}$ as \begin{align*}%\label{decomposition_local1221}
\D_{k+1|k}^{ii} &:= \nabla_{\y^i\y^i} f^i(\y^i_{k+1|k}; t_{k+1}) + \nabla_{\y^i\y^i} g^{i,i}(\y^i_{k+1|k}; t_{k+1}) \\& \hspace*{1.5cm}+ \sum_{j\in N^i} \nabla_{\y^i\y^i} g^{i,j}(\y^i_{k+1|k}, \y^j_{k+1|k}; t_{k+1})\  \\
%%\begin{equation*}%\label{decomposition_local_2290}
\B_{k+1|k}^{ij} &:= - \nabla_{\y^i\y^j}g^{i,j}(\y_{k+1|k}^i,\y_{k+1|k}^j; t_{k+1})
\end{align*}
\vspace{0.5mm}
  \STATE Compute $\displaystyle{\c_{k+1,(0)}^i \!\!=-(\D_{k+1|k}^{ii})^{-1} \nabla_{\y} F^i(\y_{k+1|k};t_{k+1})}$

    \FOR {$\tau=0,1,2,\ldots, K'-1$}
  \STATE  Exchange correction step $\c_{k,(\tau)}^i$ with neighboring nodes $j\in N^i$ 
  \STATE Compute $\c_{k+1,(\tau+1)}^i$ as % {\color{red}{??????}}
%  Execute the recursion
%  %
\begin{align*}
  \hspace{-.75cm}\c_{k+1,(\tau+1)}^i \!=\! -(\D_{k+1|k}^{ii})^{-1} 
  %\nonumber \\ &\quad  \times  
\!\Big(\!\! \sum_{j\in {N}^i } \!\! \B_{k+1|k}^{ij} \c_{k+1,(\tau)}^j \! +\! \nabla_{\y}  F^i(\y_{k+1|k};t_{k+1}) \!\Big)\nonumber
  \end{align*}
%  %
%  \vspace{-4mm}
    \ENDFOR
   \STATE Correct the trajectory prediction
  $\displaystyle{
{\y}_{k+1}^i = \y_{k+1|k}^i + \gamma \c_{k+1,(K')}^i 
}$
\end{algorithmic}
\textbf{Output}: {The corrected variable ${\y}_{k+1}^{i}$}.
\end{algorithm}
%%%%%%%%%%%%%%%%%%%%%%%%%%%%%%%%%%%
%%%%%%%%%%%%%%%%%%%%%%%%%%%%%%%%%%%
%%%%%%%%%%%%%%%%%%%%%%%%%%%%%%%%%%%

We call DPC-N as the Decentralized Prediction-Correction method that uses Newton descent in the correction step (Algorithm~\ref{alg.nc}) and the {\it exact} prediction step (Algorithm~\ref{alg.ep}) in the prediction step. We call DAPC-N as the Decentralized Approximate Prediction-Correction method that uses Newton descent in the correction step (Algorithm~\ref{alg.nc}) and the {\it approximate} prediction step (Algorithm~\ref{alg.ap}) in the prediction step. %Observe that the prediction step of the DPC-N/ DAPC-N methods can be implemented by exchanging information with neighboring nodes only as shown in Section \ref{sec:pred}. We show that the correction steps of DPC-N/DAPC-N are also implementable in a decentralized manner by invoking the result of Proposition \ref{prop.1} for the correction step.
Both DPC-N and DAPC-N require $K+K'+2$ rounds of communication per iteration.

%\subsection{Total communication counts}

For the reader's ease, we report in Table~\ref{tab.comm} the total communication counts per iteration for the presented algorithms. In particular, we report the amount of communication rounds required among the neighboring nodes, as well as the variables that have to be transmitted and the total number of scalar variables to be sent (per neighbor). 

\begin{table}
\caption{Communication requirements for the presented algorithms.}
\begin{tabular}{ccccc}
Method & & Comms. & Vars. &\hspace{-6mm} Vars. communicated \\ \toprule
\multirow{2}{*}{DPC-G/ DAPC-G} & Pred. & $K+1$ &  $\p^i_k, \y^i_k$ &\hspace{-5mm} $(K+1)p$  \\
&Corr. & $1$ & $\y^i_{k+1|k}$ &\hspace{-5mm}  $p$ \\ \toprule
\multirow{2}{*}{DPC-N/ DAPC-N} & Pred. & $K+1$ & $\p^i_k, \y^i_k$ &\hspace{-5mm} $(K+1)p$ \\ &Corr. & $K'+1$ & $\c^i_k, \y^i_{k+1|k}$ &\hspace{-5mm}  $(K'+1)p$ \\\bottomrule
\end{tabular}
\label{tab.comm}
\end{table}

\section{Convergence analysis}\label{sec:convg}

We continue by establishing the convergence of the methods presented in Section \ref{sec:algorithms}. In particular, we show that as time passes the sequence $\{\y_k\}$ approaches a neighborhood of the optimal trajectory $\y^*(t_k)$ at discrete time instances $t_k$. To establish our results, we require the following conditions.

%%%%%%%%%%%%%%%%%%%%%%%%%%%%%%%%%%%%
%%%%%%%%%%%%%%%%%%%%%%%%%%%%%%%%%%%%
%%%%%%      A  S  S  U  M  P  T  I  O  N    %%%%%%%%%%%%%%
%%%%%%%%%%%%%%%%%%%%%%%%%%%%%%%%%%%%
%%%%%%%%%%%%%%%%%%%%%%%%%%%%%%%%%%%%
%\begin{assumption}\label{as.0} 
%There exists a set $Y = Y^1 \times \dots \times Y^n \subseteq \reals^{np}$ whose interior contains the optimal argument trajectory $\y^*(t)$ of~\eqref{general_cost} for each $t$, i.e., $\y^*(t) \in \textrm{int}(Y), \,$ for $t\geq0 $. 
%\end{assumption}

%%%%%%%%%%%%%%%%%%%%%%%%%%%%%%%%%%%%
%%%%%%%%%%%%%%%%%%%%%%%%%%%%%%%%%%%%
%%%%%%      A  S  S  U  M  P  T  I  O  N    %%%%%%%%%%%%%%
%%%%%%%%%%%%%%%%%%%%%%%%%%%%%%%%%%%%
%%%%%%%%%%%%%%%%%%%%%%%%%%%%%%%%%%%%
\begin{assumption}\label{as:first} The local functions $f^i$ are twice differentiable and the eigenvalues of their Hessians $\nabla_{\y^i\y^i} f^i(\y^i;t)$ for all $i$ are contained in a compact interval $[m,M]$ with $m>0$. Hence the aggregate function $f(\y;t):=\sum_{i\in V} f^i(\y^i; t)$ has a uniformly bounded spectrum, i.e.
\begin{equation}\label{eq:eigenvalues_f}
m\mathbf{ I}\  \preceq\  \nabla_{\y\y} f(\y;t) \ \preceq\ M\mathbf{I}.
\end{equation}
\end{assumption}

%%%%%%%%%%%%%%%%%%%%%%%%%%%%%%%%%%%%
%%%%%%%%%%%%%%%%%%%%%%%%%%%%%%%%%%%%
%%%%%%      A  S  S  U  M  P  T  I  O  N    %%%%%%%%%%%%%%
%%%%%%%%%%%%%%%%%%%%%%%%%%%%%%%%%%%%
%%%%%%%%%%%%%%%%%%%%%%%%%%%%%%%%%%%%
\begin{assumption}\label{as:2}
The functions $g^{i,i}(\y^i; t)$ and $g^{i,j}(\y^i,\y^j; t)$ are twice differentiable. The Hessian of the aggregate in \eqref{cost_functions_2}, denoted as $\nabla_{\y\y} g(\y;t)$, is block diagonally dominant\cite{Feingold1962}, i.e., for all $i$,
\begin{equation}\label{eq:diag_dom}
 \left\| \nabla_{\y^i\y^i} g(\y^i, \y^j;t)^{-1} \right\|^{-1}\!\! \geq\!\!\!\! \sum_{j=1, j\neq i}^n \!\!\!\!\left\| \nabla_{\y^i\y^j} g^{i,j}(\y^i, \y^j;t)\right\| ,
\end{equation}
where by definition $\nabla_{\y^i\y^i} g(\y^i, \y^j;t) = \nabla_{\y^i\y^i} g^{i,i}(\y^i;t) + \nabla_{\y^i\y^i} g^{i,j}(\y^i, \y^j;t)$. 
The block diagonal element $\nabla_{\y^i\y^i} g(\y^i, \y^j;t)$ has eigenvalues contained in a compact interval $[\ell/2, L/2]$ with $\ell > 0$.
\end{assumption}

%%%%%%%%%%%%%%%%%%%%%%%%%%%%%%%%%%%%
%%%%%%%%%%%%%%%%%%%%%%%%%%%%%%%%%%%%
%%%%%%      A  S  S  U  M  P  T  I  O  N    %%%%%%%%%%%%%%
%%%%%%%%%%%%%%%%%%%%%%%%%%%%%%%%%%%%
%%%%%%%%%%%%%%%%%%%%%%%%%%%%%%%%%%%%
\begin{assumption}\label{as:last}
The derivatives of the global cost $F(\y;t)$ defined in \eqref{general_cost} are bounded for all $\y\in \reals^{np}$ and $ t \geq 0$ as
\begin{align}
\|\nabla_{t\y} F(\y; t)\|\!\leq \!C_0, \, \|\nabla_{\y\y\y} F(\y; t)\|\!\leq \!C_1,\, \nonumber \\  \|\nabla_{\y t\y} F(\y; t)\|\!\leq \!C_2,\, \|\nabla_{t t\y} F(\y; t)\|\!\leq \!C_3\,.  
\end{align}
%for all $\y\in Y$ and uniformly in $t$.
\end{assumption}

%%%%%%%%%%%%%%%%%%%%%%%%%%%%%%%%%%%%
%%%%%%%%%%%%%%%%%%%%%%%%%%%%%%%%%%%%
%%%%%%      M  A  I  N      M  A  T  T  E  R    %%%%%%%%%%%%%
%%%%%%%%%%%%%%%%%%%%%%%%%%%%%%%%%%%%
%%%%%%%%%%%%%%%%%%%%%%%%%%%%%%%%%%%%

%Assumption~\ref{as.0} is a weak assumption and for the case that the set $Y$ is $\mathbb{R}^{np}$, we only assume the existence of a solution for~\eqref{general_cost} at each time $t$. However, it is very useful in practice, when we know a priori that the solution trajectory has to be, for instance, positive.  
%%%
From the bounds on the eigenvalues of Hessians $\nabla_{\y\y} f(\y;t)$ and $\nabla_{\y\y} g(\y; t)$ in Assumptions~\ref{as:first} and \ref{as:2}, respectively, and from the block diagonal Gerschgorin Circle Theorem~\cite{Feingold1962} it follows that the spectrum of $\nabla_{\y\y} g(\y;t)$ lies in the compact set $[0, L]$, and the one of the Hessian of the global cost  $\nabla_{\y\y} F(\y; t)$ uniformly satisfies
%%%
\begin{equation}\label{eig_bounds}
m\ \!\mathbf{ I}\  \preceq\  \nabla_{\y\y} F(\y; t)\ \preceq\ (L+M)\ \!\mathbf{I} \; .
\end{equation}
%%%
Assumptions \ref{as:first} and \ref{as:2}, besides guaranteeing that the problem stated in~\eqref{general_cost} is strongly convex and has a \emph{unique} solution for each time instance, imply that the Hessian $ \nabla_{\y\y} F(\y; t)$ is invertible. Moreover, the higher-order derivative bounds imply the Lipschitz continuity of the gradients, Hessians, and mixed partial derivatives of the Hessians. These conditions, in addition to higher-order derivative conditions on $F$, as in Assumption~\ref{as:last}, frequently appear in the analysis of methods for time-varying optimization, and are required to establish convergence~\cite{Dontchev2013, Jakubiec2013, Ling2013}. 

Assumptions~\ref{as:first} and ~\ref{as:last} are sufficient to show that the solution \emph{mapping} $t \mapsto \y^*(t)$ is single-valued and locally Lipschitz continuous in $t$, and in particular, 
\begin{equation}\label{eq.lip}
\|\y^*(t_{k+1}) \!-\! \y^*(t_{k})\| \leq \! \frac{1}{m}\|\nabla_{t\y} F(\y; t)\| (t_{k+1}\!-\!t_{k}) \leq \frac{C_0 h }{m}, 
\end{equation} 
see for example~\cite[Theorem 2F.10]{Dontchev2009}. This gives us a link between the sampling period $h$ and the allowed variations in the optimizers.  This also gives a better understanding on the time-varying assumptions on the uniform boundedness of the time derivatives of the gradient $\nabla_{t\y}F(\y; t)$ and $\nabla_{tt\y}F(\y;t)$. In particular the bounds $C_0$ and $C_3$ require that the change and the rate of change of the optimizer be bounded. If the optimizer were the position of a moving target to be estimated, then $C_0$ and $C_3$ would be a bound on its velocity and acceleration. Finally, the bound on $\nabla_{\y t\y}F(\y; t)$ means that the quantity $\nabla_{t\y} F(\y; t)$ is Lipschitz continuous w.r.t. $\y$ uniformly in $t$. That is to say that close by points $\y$ and $\y'$ need to have similar gradient time-derivatives: e.g., if the target position is perturbed by a small amount $\delta \y$ then its velocity is perturbed by an amount not bigger than $C_2 \delta \y$.

\subsection{Discrete sampling error}

We start the convergence analysis by deriving an upper bound on the norm of the approximation error $\mathbold{\Delta}_k\in \reals^{np}$ that we estimate through a Taylor approximation in~\eqref{tay}. The error is defined as the difference between the predicted $\y_{k+1|k}$ in~\eqref{eq.qu} (with $\y_k = \y^*(t_k)$) and the exact prediction $\y^*(t_{k+1})$, starting from the same initial condition $\y^*(t_k)$, i.e.,
%%%
\begin{align}\label{eq.def.delta}
\mathbold{\Delta}_k : = \y_{k+1|k} - \y^*(t_{k+1}) .
\end{align}
In the following proposition, we upper bound the norm $\|\mathbold{\Delta}_k \|$ of the discretization error, which encodes the error due to the prediction step and is central to all our convergence results.

%%%%%%%%%%%%%%%%%%%%%%%%%%%%%%%%%%%%
%%%%%%%%%%%%%%%%%%%%%%%%%%%%%%%%%%%%
%%%%%     P  R   O   P  O  S  I  T   I  O  N     %%%%%%%%%%%%%
%%%%%%%%%%%%%%%%%%%%%%%%%%%%%%%%%%%%
%%%%%%%%%%%%%%%%%%%%%%%%%%%%%%%%%%%%
\vspace{1mm}
\begin{proposition}\label{prop.err}
Let Assumptions~\ref{as:first}-\ref{as:last} hold true. Define the discretization constant $\Delta$ as
%
%\begin{equation}\label{def_delta}
$\Delta = ({C_0^2 C_1})/{2 m^3} + ({C_0 C_2})/{m^2} +({C_3})/{2 m}.$
%\end{equation}
%
The norm of $\mathbold{\Delta}_k$ in~\eqref{eq.def.delta} is upper bounded by
\begin{equation}\label{prop_claim_err_bound}
\|\mathbold{\Delta}_k\| \leq   \Delta \ \! h^2 = O(h^2). 
\end{equation}
\end{proposition}
%
%\vspace{1mm}
%
%\begin{myproof}
% by substituting the function $f(\x; t)$ with the function $F(\y; t)$. 
%\end{my proof}
%
%\vskip1mm
%
Proposition~\ref{prop.err}, which is established as Proposition~1 in \cite{Paper1}, states that the norm of the discrete sampling error $\|\mathbold{\Delta}_k\|$ is bounded above by a constant which is in the order of $O(h^2)$. %We use this bound in establishing convergence of the proposed methods.

%\subsection{Approximate Hessian inverse error}

A second source of error to take into account, when studying the asymptotic behavior of the algorithms in Section \ref{sec:algorithms}, 
is the error due to approximating the Hessian inverse by a truncated Taylor expansion in~\eqref{eq:approx}. We bound this error as  a function of the approximation level $K$, which is the number of communication rounds among neighboring nodes.
%
%\vskip1mm
\begin{proposition}\label{th.hessian}
Under Assumptions~\ref{as:first} and \ref{as:2}, the $K$-th order approximate inverse Hessian in~\eqref{eq:approx} is well-defined. In addition, its eigenvalues are upper bounded as
\begin{equation}\label{def.H}
\|\appH\| \leq  H:= \frac{m + L/2}{m(m+\ell/2)}.
\end{equation} 
Furthermore, if we define the error of the Hessian inverse approximation as 
%\begin{equation}
$e_k = \|{\bf I} - \nabla_{\y\y}F(\y_k; t_k) {\H}_{k, (K)}^{-1}\|$,
%\end{equation}
the error $e_k$ is bounded above as
\begin{equation}\label{eq.errorE}
e_k \leq \varrho^{K+1}, \quad \textrm{where} \,\, \varrho := ({L/2})/({m + L/2}).
\end{equation}
\end{proposition}
%
%\vskip1mm
\begin{myproof}
See Appendix~\ref{app.hessian}. %\textcolor{blue}{cite Aryan's paper instead???}.
\end{myproof}
%\vskip1mm
%
%\textcolor{blue}{I think this result is exactly the same as Proposition 3 in Network Newton part 1, and therefore we can quote its proof without redoing it in Appendix \ref{app.hessian}. Let me know if you agree, and then we can cut out Appendix \ref{app.hessian}.}

Besides quantifying the error coming from approximating the Hessian inverse, Proposition \ref{th.hessian} provides trade-offs between communication cost and convergence accuracy. It shows that a larger $K$ leads to more accurate approximation of the Hessian inverse at the price of more communications. 

\subsection{Gradient tracking convergence}

In the following theorem, we establish that the sequence generated by the DPC-G and DAPC-G algorithms asymptotically converges to a neighborhood of the optimal trajectory whose radius depends on the discretization error.

%%%%%%%%%%%%%%%%%%%%%%%%%%%%%%%%%%%%
%%%%%%%%%%%%%%%%%%%%%%%%%%%%%%%%%%%%
%%%%%%      T   H   E   O   R   E   M      %%%%%%%%%%%%%%%
%%%%%%%%%%%%%%%%%%%%%%%%%%%%%%%%%%%%
%%%%%%%%%%%%%%%%%%%%%%%%%%%%%%%%%%%%
%\vskip1mm
\begin{theorem}\label{theorem.gradient}
Consider the sequence $\{\y_k\}$ generated by the DPC-G or DAPC-G algorithm, which uses Algorithm~\ref{alg.ep} (or \ref{alg.ap}) as prediction step and Algorithm~\ref{alg.gc} as correction step. Let Assumptions~\ref{as:first}-\ref{as:last} hold and define constants $\rho$ and $\sigma$ as
\begin{align}\label{def_sigma}
\rho\! :=\! \max\{|1\!-\!\gamma m|,|1-\gamma (L\!+\!M)|\}, 
%\nonumber \\ 
\sigma\! :=\! 1\! +\! h\left[\!\frac{C_0 C_1}{m^2}\! + \!\frac{C_2}{m}\!\right].
\end{align}
%%
%Define the boolean variable $\booa$ to be $1$ if we run the DAPC-G algorithm and $0$ otherwise. 
Further, recall the definition of $\varrho$ in \eqref{eq.errorE} and define the function $\Gamma: (0,1)\times\mathbb{N} \to \reals$ as
$\Gamma(\varrho, K) = ({C_0}/{m}) \varrho^{K+1}$. Choose the step-size as 
\begin{equation}
\gamma < 2/(L+M),
\end{equation}
so that $\rho < 1$. Then,
\begin{enumerate}[i)]
%%%%
\item For any sampling period $h$, the sequence $\{\y_k\}$ converges to $\y^*(t_k)$ Q-linearly up to a bounded error, as
%%%%%
\begin{align}\label{result1}
\hskip-1cm \limsup_{k\to \infty}\|{\y}_{k} - \y^*(t_{k})\| &= O(h) + O(h \Gamma(\varrho, K)) + O(h^2).
\end{align}
%%%%
\item If the sampling period $h$ is chosen such that
\begin{equation}\label{eq.hh}
h < \left[\frac{C_0 C_1}{m^2} + \frac{C_2}{m}\right]^{-1} (\rho^{-1} - 1),
\end{equation}
then the sequence $\{\y_k\}$ converges to $\y^*(t_k)$ Q-linearly up to a bounded error as
%%%%
\begin{align}\label{result2}
\hskip-1cm \limsup_{k\to \infty}\|{\y}_{k} - \y^*(t_{k})\| &= O(h \Gamma(\varrho, K)) + O(h^2).
\end{align}
%%%%
\end{enumerate} 
\end{theorem}
\begin{myproof}
See Appendices~\ref{app.DPC-Gconvg}-\ref{app.DAPC-Gconvg}, where the error bounds and convergence rate constant are explicitly computed in terms of the functional bounds of Assumptions~\ref{as:first}-\ref{as:last}.
\end{myproof} 

Theorem \ref{theorem.gradient} establishes the convergence properties of DPC-G and DAPC-G for particular parameter choices. In both cases, the linear convergence to a neighborhood is shown, provided the step-size satisfies $\gamma<2/(L+M)$. Moreover, the accuracy of convergence depends on the choice of the sampling period $h$, and for any sampling period, the result in \eqref{result1} holds. In this case the accuracy of convergence is of the order $O(h)$. If the sampling period $h$ is chosen such that  $\rho \sigma < 1$ (that is~\eqref{eq.hh} holds), then the result in \eqref{result2} is valid. 

If $\rho \sigma < 1$ is satisfied and %$\Gamma(\varrho, K) = O(h)$ holds, i.e, the %approximation level of the Hessian inverse is chosen as
%\begin{equation}\label{Kcond}
%$K \gg \lceil \log h/ \log \varrho -1 \rceil,$
%\end{equation}
%then the error bound in~\eqref{result2} becomes
%%%%%
%\begin{align}\label{asy_3}
%\limsup_{k\to\infty} \|\y(t_k) - \y^*(t_k)\| \leq O(h^2).
%\end{align}
%%%%%
%I.e., if 
the approximation level $K$ %(or equivalently the number of communication rounds among neighboring nodes) 
is chosen sufficiently large, then $\Gamma(\varrho, K)$ is negligible and we regain an error bound of $O(h^2)$, which is compatible with centralized algorithms \cite{Paper1}.% we prove that the centralized versions of DPC-G and DAPC-G reach at most a $O(h^2)$ asymptotic error bound.

\subsection{Newton tracking convergence}
We turn to analyzing the DPC-N and DAPC-N algorithms.%decentralized trajectory tracking methods which use second-order information in the correction step. %In particular, we next state the convergence guarantees of the DPC-N algorithm. 
%%%%%%%%%%%%%%%%%%%%%%%%%%%%%%%%%%%%%
%%%%%%%%%%%%%%%%%%%%%%%%%%%%%%%%%%%%%
%%%%%%    T  H  E  O   R   E  M     %%%%%%%%%%%%%%%%
%%%%%%%%%%%%%%%%%%%%%%%%%%%%%%%%%%%%%
%%%%%%%%%%%%%%%%%%%%%%%%%%%%%%%%%%%%%
%\vskip1mm 
\begin{theorem}\label{th.newton}
Denote $\{\y_k\}$ as the sequence generated by the DPC-N or DAPC-N method, which respectively uses Algorithm~\ref{alg.ep} or~\ref{alg.ap} as its prediction, and Algorithm~\ref{alg.nc} as its correction. Let Assumptions~\ref{as:first}-\ref{as:last} hold and fix $K$ and $K'$ as the Hessian inverse approximation levels for the prediction and correction steps, respectively, with the function $\Gamma$ defined as in Theorem \ref{theorem.gradient}. Fix the step-size as $\gamma\in(0,1]$.
%
%Let $\tau<1$ be a positive scalar and fix the levels of approximation of the Hessian inverse for the prediction and correction steps as $K$, and $K'$, respectively. Define constants $\alpha_1$ and $D$ as
%
%\begin{align}\label{def_alpha}
%\alpha_1 &:= \sigma\left[\frac{C_1 h D}{m}  + \frac{L+M}{C_0}\Gamma(\varrho, K')\right], \\
%D &:= \Gamma(\varrho, K) + \Delta h ,
%\end{align}
%
%where $\Delta$, $\sigma$, $\Gamma$, and $\varrho$ are defined in~Proposition \ref{prop.err}, \eqref{def_sigma}, and Theorem \ref{theorem.gradient}, respectively. 
%
There exist bounds $\bar{K}$, $\bar{h}$, and $\bar{R}$, such that if the sampling rate $h$ is chosen as $h \leq \bar{h}$, $K$ and $K'$ are chosen as $K,K'\geq \bar{K}$, and the initial optimality gap satisfies $\|\y_0-\y^*(t_0)\|\leq \bar{R}$,
%
%\begin{equation}\label{conv_region}
%\|\y_0-\y^*(t_0)\| \leq \frac{2 m}{\sigma^2 C_1}(\tau - \alpha_1) \; ,
%\end{equation} 
%
%and the sampling period $h$, and Hessian approximation levels $K$, $K'$ are chosen such that
%
%\begin{equation}\label{alpha_cond}
%$\alpha_1 < \tau $,
%\end{equation}
%
then $\{\y_k\}$ converges Q-linearly to the solution trajectory $\y^*(t_k)$ up to a bounded error as
\begin{align}\label{result_37}
\!\limsup_{k\to\infty}  \|{\y}_{k}\! -\! \y^*(t_{k})\| \!&=\! O(h \Gamma(\varrho, K) [\gamma\Gamma(\varrho, K') \!+\! 1\!-\!\gamma]) \nonumber
\\&\!\!\!\!\!\!+ O(h^2 [\gamma\Gamma(\varrho, K')\! +\! \gamma\Gamma(\varrho, K)^2\! 
+\! 1\!-\!\gamma]) \nonumber
\\&\!\!\!\!\!\!+ O(h^3 \gamma \Gamma(\varrho, K)) + O(h^4\,\gamma).  
\end{align}
In addition, if the step-size $\gamma$ is chosen arbitrarily small, the attraction region $\bar{R}$ can be made arbitrarily large. 

%
%{\color{red}provided that the initial optimality gap $\|\y_0-\y^*(t_0)\|$ is sufficiently small.}
\end{theorem}

\begin{myproof}
See Appendices~\ref{newton.proof1}-\ref{newton.proof2}. The proof is constructive, thus we also characterize the bounds on the sampling period, approximation levels, the attraction region, and finally, the constants in the asymptotic error and in the linear convergence rate. 
\end{myproof}

%\vskip1mm

Theorem~\ref{th.newton} states that 
DPC-N/DAPC-N converge to a bounded tracking error defined in~\eqref{result_37} once the algorithm reaches an attractor region. The error bound in~\eqref{result_37} depends, as expected, on the sampling period $h$ and the approximation levels $K$ and $K'$. In the worst case, the asymptotic error floor will be of the order $O(h)$. However, in some cases we may achieve tighter tracking guarantees. For example, if the approximation level $K$ and $K'$ are chosen sufficiently large, then the terms $\Gamma(\rho, K)$ and $\Gamma(\rho, K')$ are negligible, yielding
\begin{equation}\label{eq:DPC-N_asymptotic2}
\limsup_{k\rightarrow \infty} \|{\y}_{k} - \y^*(t_{k})\| = O(h^2[1-\gamma]) + O(h^4\,\gamma) \; .
\end{equation}
This is to say that the asymptotic dependence of the error of DPC-N/DAPC-N on the sampling period $h$ varies from a worst-case $O(h)$ to as tight as $O(h^4)$ (for the selection $\gamma =1$).

\begin{remark}\emph{(Step-size choice)}
The step-size choice of the presented Newton correction methods affects the convergence attraction region and the convergence speed. For large enough $K, K'$ and small enough $h$, if we choose $\gamma =1$, we obtain a standard Newton method with convergence region $\bar{R} =  2 m/ C_1 \sigma^2$, the fastest convergence speed, and smallest asymptotical error $O(h^4)$. This convergence region is larger depending on how small $C_1$ is: for quadratic functions $C_1 = 0$, and the convergence is global. If we choose $\gamma \ll 1$, then the convergence region is $\bar{R} =  2 m (\tau-1+\gamma)/ \gamma C_1 \sigma^2$, where $\tau$, with $1-\gamma<\tau<1$, is the linear convergence rate [Cf. Appendix~\ref{newton.proof1}]. This means that the attraction region can be made arbitrarily big, while the convergence rate is made smaller and smaller, and the asymptotical error is $O(h^2)$. Finally, an interesting choice is $\gamma = h \leq 1$: when $h$ is sufficiently small, then the Newton prediction/correction approximate well a continuous-time algorithm and the convergence (albeit made slow) is global. 
In practice, the choice of the step-size depends on the application at hand. One can even decide to run the DPC-G algorithm till convergence and then switch to DPC-N as an hybrid scheme, or to adopt an increasing step-size selection. These extensions are left as future research.%, while in the simulation section, we analyze the choice $\gamma =1$. 
\end{remark}

\section{Numerical evaluation}\label{sec:num}

We turn to studying the empirical validity of the performance guarantees established in Section \ref{sec:convg}. In particular, we consider the resource allocation problem in a network of interconnected devices, as in Example \ref{example2} of Section \ref{sec:prob}. %Here each sensor is a node $i\in V$ associated with a graph $\ccalG=(V,E)$, and aims to collaboratively solve a time-varying resource allocation problem. 
As presented in~\eqref{res_all_1}, the local objective functions $f^i(\y^i; t)$ represent a time-varying utility indicating the quality of transmission at a particular device $i$ and the constraints represent channel rate and capacity constraints. These constraints depend on the connectivity of the network, which is encoded in the augmented incidence matrix $\A$, defined following \eqref{res_all_1}.

By adopting an approximate augmented Lagrangian method, we obtain \eqref{res_all_2} which is an instance of~\eqref{general_cost}. Consider the case where decisions are the variables $\y^i \in \reals^{p}$, $p = 10$, for which each local utility $f^i(\y^i; t)$ associated with sensor $i$ is given as
\begin{multline}\label{eq:example1_special_case}
f^i(\y^i; t) = \frac{1}{2} (\y^i - \c^i(t))^\transp \Q^i (\y^i - \c^i(t)) + \\ \sum_{l=1}^p \log\left[1 + \exp\left(b^{i,l}(y^{i,l} - d^{i,l}(t))\right) \right],
\end{multline}
where $\y^{i,l}$ indicates the $l$-th component of the $i$-th decision variable $\y^i$, while $\Q^i \in \reals^{p\times p}, b^{i,l}\in \reals, \c^i(t)\in \reals^{p}, d^{i,l}(t)\in \reals$ are (time-varying) parameters. Straightforward computations reveal that the second order derivative of $f^i$ with respect to $\y^i$ is 
%
%\begin{equation}
%\frac{\partial^2}{\partial (y^i)^2} f^i(y^i; t) = a^i  + (b^i)^2 \frac{\left[\exp\left(b^i(y^i - d^i(t))\right) \right]}{\left[1 + \exp\left(b^i(y^i - d^i(t))\right) \right]^2},
%\end{equation}
%and it is 
contained in the bounded interval $[\lambda_{\min}(\Q^i),\lambda_{\max}(\Q^i) + \max_{l}\{(b^{i,l})^2/4\}]$.

Experimentally, we consider cases were each $\Q^i$ and $b^{i,l}$ are selected uniformly at random, and in particular $\Q_i = \textrm{diag}(\mathcal{U}^p_{[1,2]}) + \v^i {\v^i}^\transp$, with $\v_i \sim \mathcal{N}^p_{0,1}$ (that is $\v_i$ is a random vector drawn from a Gaussian distribution of mean zero and standard deviation one). With this choice $\Q^i$ is positive definite. In addition $b^{i,j} \sim \mathcal{U}^1_{[-2,2]}$. Finally, $\c^i(t)$ and $d^{i,l}(t)$ are the time-varying functions
\begin{subequations}
\begin{align}
c^{i,l}(t) = 10 \cos(\theta^{i,l}_c + \omega\,t), \quad \theta^{i,l}_c \sim \mathcal{U}[0, 2\pi),\\
d^{i,l}(t) = 10 \cos(\theta^{i,l}_d + \omega\,t), \quad \theta^{i,l}_d \sim \mathcal{U}[0, 2\pi),
\end{align} 
\end{subequations}
with $\omega = 0.1$. The sensors in the $n=50$ node wireless network are deployed randomly in the area $[-1,1]^2$ and can communicate if they are closer than a range of $r = 2.5\sqrt{2}/\sqrt{n}$, which generates a network of $l$ links. We set the vector of rate and capacity constraints to $\b = {\bf 0}$ yielding a dynamic network flow problem, with approximation level $\beta = \sqrt{20}$. 
%\textcolor{blue}{How is the region $\ccalA$ defined in which the sensors are deployed? Then, what is the strategy for spatially deploying the nodes? Once these two things are defined, defining a distance-based connectivity rule makes sense. What is $\ell$ here, and how is it different from $m$, the number of edges in $\ccalG$?  }
%We can compute the bounds in the Assumptions~\ref{as.smooth2}-\ref{as.smooth4} as
%\begin{equation}
%m = 1, \, M = 2, \, L = 1.4, \, C_0 = 2, \, C_1 = 0.8, \, C_2 = 0.8.
%\end{equation}

\subsection{Comparisons in absolute terms}
We first analyze the behavior of DPC-G, DAPC-G, DPC-N, and DAPC-N with respect to the decentralized running gradient method of~\cite{Simonetto2014c}. %\footnote{We do not consider the works in \cite{Jakubiec2013} and \cite{Ling2013}, since we employ an approximate augmented Lagrangian framework to reduce the computation/communication requirements and it would not make much sense to use these computationally more intensive algorithms. Please refer to~\cite{Paper1Camsap} for a different numerical example and comparison with these methods. {\color{red}{Do we need this footnote?}}}
%
 %However, one may derive the fact that these methods achieve an $O(h)$ tracking error, and thus behave comparably to running gradient methods.}.
%
%
Unless otherwise stated the DPC-N and DAPC-N algorithms run with unitary step-size ($\gamma = 1$).

%%%%%%%%%%%%%%%%%%%%%%%%%%%%%%%%%%%%%%%%%%%%%%%%%%%%%%%%%%%%%%%%%%%%%%%%%%%%%%%%%%%%%%%%%%%%%%%%%%%%
%%%%%%%%%%%%%%%%%%%% FIGURE
%%%%%%%%%%%%%%%%%%%%%%%%%%%%%%%%%%%%%%%%%%%%%%%%%%%%%%%%%%%%%%%%%%%%%%%%%%%%%%%%%%%%%%%%%%%%%%%%%%%%
\begin{figure}
%\psfragfig*[width=.5\textwidth, height = 5.6cm]{Figure1rev}
%{
\scriptsize
\psfrag{RunningGradientGradient}{Running gradient}
\psfrag{DeGT1}{DPC-G, $K$=$3$}
\psfrag{DeGT2}{DPC-G, $K$=$5$}
\psfrag{DeAGT1}{DAPC-G, $K$=$5$}
\psfrag{DeNT1}{DPC-N, $K$=$K'$=$3$}
\psfrag{DeNT2}{DPC-N, $K$=$K'$=$5$}
\psfrag{DeANT1}{DAPC-N, \!$K$=$K'\!$=$5$}
\psfrag{y}[c]{Error\, $\|\y_k - \y^*(t_k)\|$}
\psfrag{x}[c]{Sampling time instance $k$}
%}
\includegraphics[width=.5\textwidth, height = 5.3cm]{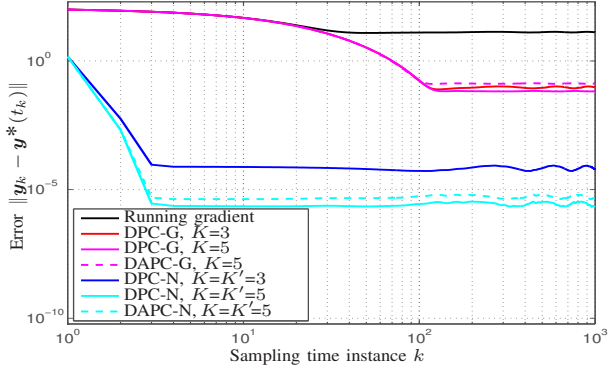}
\vskip-2mm
\caption{Error with respect to the sampling time instance $k$ for different algorithms applied
to the continuous-time sensor network resource allocation problem~\eqref{res_all_2}, with sampling interval $h = 0.1$. %DPC-N and DAPC-N achieve superior performance for $K=K'=5$.%
}
\label{fig:1}
\end{figure}

%%%%%%%%%%%%%%%%%%%%%%%%%%%%%%%%%%%%%%%%%%%%%%%%%%%%%%%%%%%%%%%%%%%%%%%%%%%%%%%%%%%%%%%%%%%%%%%%%%%%
%%%%%%%%%%%%%%%%%%%%%%%%%%%%%%%%%%%%%%%%%%%%%%%%%%%%%%%%%%%%%%%%%%%%%%%%%%%%%%%%%%%%%%%%%%%%%%%%%%%%
%%%%%%%%%%%%%%%%%%%%%%%%%%%%%%%%%%%%%%%%%%%%%%%%%%%%%%%%%%%%%%%%%%%%%%%%%%%%%%%%%%%%%%%%%%%%%%%%%%%%

%%%%%%%%%%%%%%%%%%%%%%%%%%%%%%%%%%%%%%%%%%%%%%%%%%%%%%%%
%%%%%%%%%%%%%%%%%%%% FIGURE
%%%%%%%%%%%%%%%%%%%%%%%%%%%%%%%%%%%%%%%%%%%%%%%%%%%%%%%%%%%%%%%%%%%%%%%%%%%%%%%%%%%%%%%%%%%%%%%%%%%%%%%%%%%%%%%%%%%%%%%%%%%%%%%%%%%%%%%%%%%%%%%%
\begin{figure}[t]
%\psfragfig*[width=.5\textwidth, height = 5.6cm]{Figure2rev}
%{
\scriptsize
\psfrag{RunningGradientGradient}{Running gradient}
\psfrag{DeGT1}{DPC-G, $K$=$3$}
\psfrag{DeGT2}{DPC-G, $K$=$5$}
\psfrag{DeAGT}{DAPC-G, $K$=$5$}
\psfrag{GTT}{DPC-G, $K\to\infty$}
\psfrag{DeNT1}{DPC-N, $K$=$K'$=$3$}
\psfrag{DeNT2}{DPC-N, $K$=$K'$=$5$}
\psfrag{DeANT}{DAPC-N, \!\!$K$=$K'\!$=$5$}
\psfrag{NTT}{DPC-N, \!\!$K$=$K'\!\!\!\to\!\!\infty$}
\psfrag{O1}{$O(h)$}\psfrag{O2}{$O(h^2)$}\psfrag{O4}{$O(h^4)$}
\psfrag{y}[c]{Asymptotic error bound}
\psfrag{x}[c]{Sampling period $h$}
%}
\includegraphics[width=.5\textwidth, height = 5.3cm]{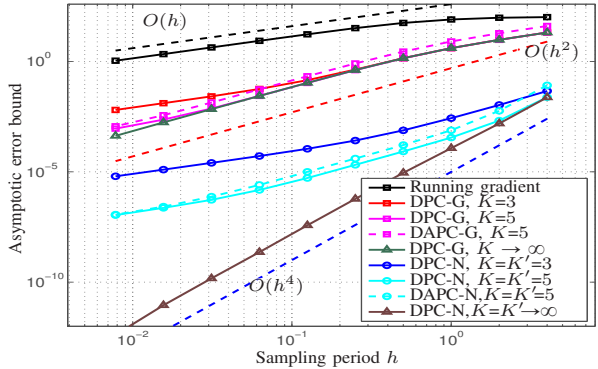}
\vskip-2mm
\caption{Asymptotic error bound $\max_{k>\bar{k}} \left\{\|\y_k - \y^*(t_k)\|\right\}$ as compared with the sampling interval $h$. Dotted straight lines represent error bounds $O(h^r)$ for $r=1,2,4$.}
%
%for $\bar{k}=800$ for $h\geq 1/16$ or $\bar{k}=2000$ for $h<1/16$.
%
% We see that DPC-N and DAPC-G achieve an error bound orders of magnitude smaller than competing methods when applied to the problem~\eqref{res_all_2}.}
\label{fig:2}
\end{figure}

In Figure~\ref{fig:1}, we depict how the different algorithms reach convergence as time passes for a fixed sampling interval of $h = 0.1$. Observe that the running gradient method achieves the worst tracking performance of around $\|\y_k - \y^*(t_k)\|\approx 10$, whereas DPC-G for various levels of communication rounds in the prediction and correction steps $K$ and $K'$ achieves an error near $10^{-1}$. Using second-order information in the correction step, as with DPC-N and DAPC-N, achieves superior performance, with tracking errors of at least $\|\y_k - \y^*(t_k)\|\approx 10^{-5}$. Moreover, the time-approximation in DAPC-G and DAPC-N does not degrade significantly the asymptotic error, while the number of communication rounds $K$ and $K'$ play a more dominant role, especially in the case of DPC-N. 

We also observe this trend in Figure~\ref{fig:2}, where we analyze the behavior varying the sampling period $h$. We approximate the asymptotic error bound as 
$%\begin{equation}
\max_{k>\bar{k}} \left\{\|\y_k - \y^*(t_k)\|\right\}
$, %\end{equation}
  for a given $\bar{k}$, where we set $\bar{k}=800$ for $h\geq 1/16$ or $\bar{k}=2000$ for $h<1/16$. We may observe empirical confirmation of the error bounds established by Theorems \ref{theorem.gradient} and \ref{th.newton} in Section~\ref{sec:convg}. In particular, the running gradient has an asymptotic error approximately as $O(h)$, whereas that of DPC-G varies between $O(h)$ and $O(h^2)$ depending on the approximation level $K$ and $h$. Moreover, DPC-N achieves an asymptotic tracking error varying between $O(h)$ and $O(h^4)$. 
  
\begin{figure}[t]
%\psfragfig*[width=.5\textwidth, height = 5.6cm]{Figure3rev}
%{
\scriptsize
\psfrag{RunningGradientRunningGradient}{DPC-G, $K$=$3$}
\psfrag{G1}{DPC-N, $K$=$K'$=$3$, $\gamma$=$.1$}
\psfrag{G2}{DPC-N, $K$=$K'$=$3$, $\gamma$=$.5$}
\psfrag{G3}{DPC-N, $K$=$K'$=$3$, $\gamma$=$1$-$.9/k$}
\psfrag{G4}{DPC-N, $K$=$K'$=$3$, $\gamma$=$1$}
\psfrag{y}[c]{Error\, $\|\y_k - \y^*(t_k)\|$}
\psfrag{x}[c]{Sampling time instance $k$}
%}
\includegraphics[width=.5\textwidth, height = 5.3cm]{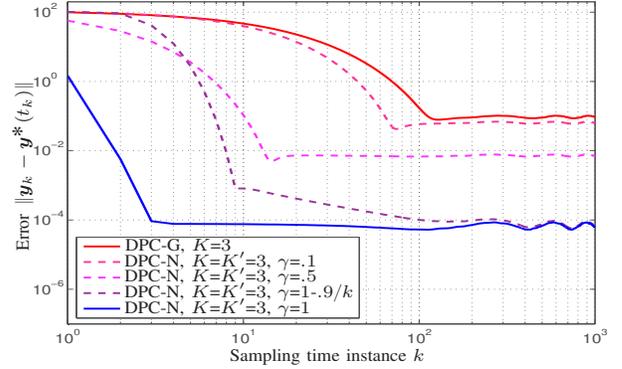}
\vskip-2mm
\caption{Error with respect to the sampling time instance $k$ for DPC-G and DPC-N with different step-size $\gamma$, applied
to the continuous-time sensor network resource allocation problem~\eqref{res_all_2}, with sampling interval $h = 0.1$. \color{red} }
\label{fig:3}
\end{figure}

In Figure~\ref{fig:3}, we depict the behavior in time for different choices of step-size $\gamma$ for DPC-N. As we notice, varying from a small step-size $\gamma = 0.1 = h$ to the biggest one of $\gamma = 1$, the convergence becomes faster (yet theoretically more local). An increasing choice of step-size as $\gamma = 1 - .9/k$ seems to combine both larger convergence region, reasonably fast convergence, and small asymptotical error.

Since DPC-N is a computationally more demanding method in terms of communication requirements and computational latency, we study the effect of fixing the former parameter.

%%%%%%%%%%%%%%%%%%%%%%%%%%%%%%%%%%%%%%%%%%%%%%%%%%%%%%%%%%%%%%%%%%%%%%%%%%%%%%%%%%%%%%%%%%%%%%%%%%%%
%%%%%%%%%%%%%%%%%%%% FIGURE
%%%%%%%%%%%%%%%%%%%%%%%%%%%%%%%%%%%%%%%%%%%%%%%%%%%%%%%%%%%%%%%%%%%%%%%%%%%%%%%%%%%%%%%%%%%%%%%%%%%%

%%%%%%%%%%%%%%%%%%%%%%%%%%%%%%%%%%%%%%%%%%%%%%%%%%%%%%%%%%%%%%%%%%%%%%%%%%%%%%%%%%%%%%%%%%%%%%%%%%%
%%%%%%%%%%%%%%%%%%%%%%%%%%%%%%%%%%%%%%%%%%%%%%%%%%%%%%%%%%%%%%%%%%%%%%%%%%%%%%%%%%%%%%%%%%%%%%%%%%

\subsection{Comparisons with fixed communication effort}

In practice the communication and computation requirements for each of the nodes of the network will be fixed by hardware and bandwidth constraints. Let us fix the time, as a percentage of the sampling period $h$, for the prediction and correction step. Let us say that we have at most a time of $r h$ ($r\leq .5$) to do prediction and $r h$ to do correction. 

Each time a new function is sampled, each of the proposed algorithms will perform a number of correction steps $n_{\textrm{C}} \geq 1$. Each of them will consist of either $n_{\textrm{C}}$ gradient steps, involving each broadcasting $p$ scalar values to the neighbors and receiving $p N_i$ scalar values from them, or $n_{\textrm{C}}$ approximate Newton steps, involving each broadcasting $p(K'+1)$ scalar values to the neighbors and receiving $p N_i (K'+1)$ scalar values from them. 

Once the corrected variable is derived, it can be implemented (e.g., generating the control action). In the remaining time, while waiting for another sampled cost function, the proposed algorithms can perform a prediction step, involving for each node broadcasting $p(K+1)$ scalar values to the neighbors and receiving $p N_i(K+1)$ scalar values from them. For the running schemes, there is no prediction, but we assume here that the variables are further optimized by other extra correction steps, and hence start at the next time with a better initialization. These further correction steps, say $n_{\textrm{EC}}$, can be gradient or Newton. 

Define $\bar{t}$ to be the time required for one round of broadcasting and receiving data from and to the neighbors and assume that it is the same for each node and it scales linearly with the number of communication rounds $K$ and $K'$ (since it has to be done sequentially). Suppose, as empirically observed, that the computation time for the nodes is negligible w.r.t. the communication time. In this context the number of correction and prediction rounds can be chosen according to the constraints on time:
\begin{subequations}\label{kk}
\begin{align}
% &&\textrm{Correction} & \textrm{Prediction or Extra correction} \nonumber \\
\textrm{(RG)} && n_{\textrm{C}} \bar{t} = r h\,  & \quad n_{\textrm{EC}}\bar{t} = r h, \\
\textrm{(RN)} && n_{\textrm{C}} (K'+1)\bar{t} = r h\,  & \quad n_{\textrm{EC}}(K+1)\bar{t} = r h, \\
\textrm{(DPC-G)} && n_{\textrm{C}} \bar{t} = r h\,  & \quad (K+1)\bar{t} = r h, \\
\textrm{(DPC-N)} && n_{\textrm{C}} (K'+1)\bar{t} = r h\,  & \quad (K+1)\bar{t} = r h,
\end{align}
\end{subequations}
where RG indicates the running gradient method and RN the running Newton. In the following simulation, we fix $r = 0.5$, $K = K'$ and $n_{\textrm{C}}$, $n_{\textrm{EC}}$ to be $1$ for RN and DPC-N. We fix $\bar{t} = 1/10$~s (for bigger values of $p$, that is the dimension of the decision variable, this time will be longer, and vice-versa).

\begin{figure}
%\psfragfig*[width=.5\textwidth, height = 5.6cm]{Figure4rev}
%{
\scriptsize
\psfrag{RunningGradientGradient}{Running gradient}
\psfrag{DeGT2}{DPC-G}
\psfrag{DeGT1}{Running Newton}
\psfrag{DeAGT}{DPC-N, $\gamma$=$1$}
\psfrag{y}[c]{Asymptotic error bound}
\psfrag{x}[c]{Sampling period $h$}
%}
\includegraphics[width=.5\textwidth, height = 5.2cm]{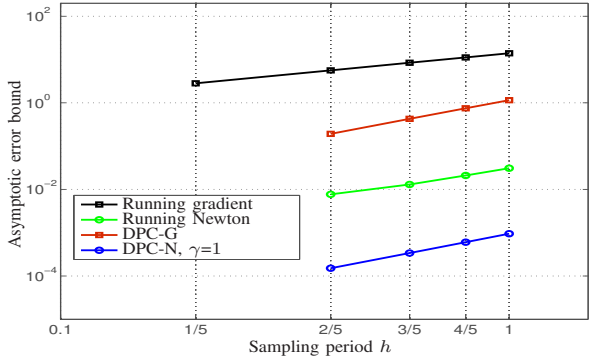}
\vskip-2mm
\caption{Asymptotic error with respect to the sampling period $h$ for different algorithms when the number of communication rounds is chosen according to bandwidth constraints as in~\eqref{kk}.}
%\textcolor{blue}{we need to fix the legend of this figure}}
\label{fig:4}
\end{figure}

In Figure~\ref{fig:4}, we report the results in terms of asymptotical error when optimizing the number of communication rounds according to~\eqref{kk}, for different sampling periods. In this context, e.g., for $h = 1$~s, we can run RG with $n_\textrm{C} = n_\textrm{EC} = 5$ correction and extra correction rounds, RN with $K = K' = 4$ communication rounds for correction and extra correction, DPC-G with $K = 4$ communication rounds for prediction and $n_\textrm{C} = 5$ rounds of correction, and DPC-N with $K = K' = 4$ communication rounds for prediction and correction, respectively. For the other values of sampling period $h$, similar calculations give us the optimized values for $n_{\textrm{C}}$, $n_{\textrm{EC}}$, $K$, and $K'$.  
As we observe, when the sampling period is big enough, so that DPC-N is implementable, then it seems to be the best strategy to go for. We also notice that running gradient, even with the extra correction steps, should be avoided unless all the remaining algorithms are unviable (as for $h = 1/5$~s).

%%%%%%%%%%%%%%%%%%%%%%%%%%%%%%%%%%%%%%%%%%%%%%%%%%%%%%%%%%%%%%%%%%%%%%%%%%%%%%%%%%%%%%%%%%%%%%%%%%%
%%%%%%%%%%%%%%%%%%% FIGURE
%%%%%%%%%%%%%%%%%%%%%%%%%%%%%%%%%%%%%%%%%%%%%%%%%%%%%%%%%%%%%%%%%%%%%%%%%%%%%%%%%%%%%%%%%%%%%%%%%%

% !TEX root = TimeVarying_part2.tex
\section{Conclusions}\label{sec:conclusion}

We considered continuously varying convex programs whose objectives may be decomposed into two parts:  a sum of locally available functions at the nodes and a part that is shared between neighboring nodes. To solve this problem and track the solution trajectory, we proposed a decentralized iterative procedure which samples the problem at discrete times. Each node predicts where the solution trajectory will be at the next time via an approximation procedure in which it communicates with its neighbors, and then corrects this prediction by incorporating information about how the local objective is varying, again via a decentralized local approximation. We developed an extension of this tool which allows for the case when the dynamical behavior of the objective must be estimated.

 \begin{table}\centering\hfill
\caption{Summary of proposed methods and convergence results.}
\label{tab.1}
\begin{tabular}{ccccc}
Method & DPC-G& DAPC-G& DPC-N& DAPC-N\\
\toprule 
Prediction & Alg.~\ref{alg.ep}  & Alg.~\ref{alg.ap} & Alg.~\ref{alg.ep}& Alg.~\ref{alg.ap} \\ 
Correction & Alg.~\ref{alg.gc} &  Alg.~\ref{alg.gc} & Alg.~\ref{alg.nc} & Alg.~\ref{alg.nc} \\ 
Best Error bound & $O(h^2)$  & $O(h^2)$ &  $O(h^4)$   & $O(h^4)$  \\
%Reference & Eqns.~\eqref{result1}-\eqref{result2}  & Eqn.~\eqref{result11}-\eqref{result22}&  Eqn.~\eqref{result_37}]  & Eqn.~\eqref{result_371} \\
\bottomrule
\end{tabular}
\end{table}

We established that this decentralized approximate second-order procedure converges to an asymptotic error bound which depends on the length of the sampling interval and the amount of communications in the network. Moreover, we established that this convergence result also applies to the case where time derivatives must be approximated. A summary of the proposed methods and their performance guarantees is given in Table \ref{tab.1}. Finally, we applied the developed tools to a resource allocation problem in a wireless network, demonstrating its practical utility and its ability to outperform existing running methods by orders of magnitude.

\appendices

\section{Proof of Proposition~\ref{th.hessian} }\label{app.hessian}

We generalize the proofs of Propositions 2 and 3 in~\cite{mokhtari2015network1} to establish the result. Start by defining the matrix $\hat{\D}_k = \textrm{diag}[\nabla_{\y\y}g(\y_k;t_k)]$ which is positive definite due to Assumption~\ref{as:first}. Thus, we can write
\begin{align}\label{eq:dummy0}
\!\!\!\! \|\D_k^{-1/2} \B_k \D_k^{-1/2}\|\!&= \!\!
\| \D_k^{-1/2}\hat{\D}_k^{1/2} \hat{\D}_k^{-1/2} \B_k \hat{\D}_k^{-1/2} \nonumber%\times \\ 
\hat{\D}_k^{1/2}\D_k^{-1/2} \|  \\
&\leq \|\D_k^{-1/2}\hat{\D}_k^{1/2} \|^2 \|\hat{\D}_k^{-1/2} \B_k \hat{\D}_k^{-1/2}\|.
\end{align}
where the inequality in \eqref{eq:dummy0} is implied by the Cauchy-Schwartz inequality. We proceed to bound both terms on the right-hand side of \eqref{eq:dummy0}, starting with the rightmost term. The matrix $\hat{\D}_k^{-1/2} \B_k \hat{\D}_k^{-1/2}$ is conjugate to the matrix $ \hat{\D}_k^{-1}\B_k  $, which means that the latter has the same eigenvalues of the former. By construction, the matrix $ \B_k \hat{\D}_k^{-1} $ is equivalent to 
\begin{equation}\label{eq.db}
\hat{\D}_k^{-1}\B_k  = {\bf I} - \hat{\D}_k^{-1}\nabla_{\y\y}g(\y_k;t_k).
\end{equation}
The matrix $\nabla_{\y\y}g(\y_k;t_k)$ is block diagonally dominant (Assumption~\ref{as:2}); by the definition of the matrix $\hat{\D}_k$, this means
\begin{equation}\label{eq.dd1}
\|[\hat{\D}_k^{-1}]^{ii}\|^{-1}  \geq \sum_{j=1, j\neq i}^n \!\!\left\| \nabla_{\y^i\y^j} g^{i,j}(\y^i_k, \y^j_k;t)\right\|, \, \textrm{for all} \, i.
\end{equation}
%or equivalently, 
%\begin{equation}
%1  \geq \|[\hat{\D}_k^{-1}]^{ii}\|\sum_{j=1, j\neq i}^n \!\!\!\!\left\| \nabla_{\y^i\y^j} g^{i,j}(\y^i_k, %\y^j_k;t)\right\|.
%\end{equation}
Now consider the matrix $\hat{\D}_k^{-1}\nabla_{\y\y}g(\y_k;t_k)$, by the block Gershgorin Circle theorem~\cite{Feingold1962} its eigenvalues are contained in the circles defined by all the $\mu$'s that verify
\begin{equation}
 \!\!\|({\bf I} - \mu {\bf I})^{-1}\|^{-1}   \!\!\leq \!\!\!\!\sum_{j=1, j\neq i}^n \!\!\!\!\left\| [\hat{\D}_k^{-1}]^{ii} \nabla_{\y^i\y^j} g^{i,j}(\y^i_k, \y^j_k;t)\right\|  \!\!\leq  \!\!1,
\end{equation}
where the last inequality comes from~\eqref{eq.dd1}. Therefore the eigenvalues are contained in the compact set $[0,2]$. This means that the matrix $\hat{\D}_k^{-1}\B_k$ in~\eqref{eq.db} has eigenvalues contained in the compact set $[-1,1]$. Taken with the fact that the Frobenius norm of the matrix $\hat{\D}_k^{-1/2} \B_k \hat{\D}_k^{-1/2}$ is bounded above by its maximum eigenvalue, we have
\begin{equation}\label{eq:dummy21}
\|\hat{\D}_k^{-1/2} \B_k \hat{\D}_k^{-1/2}\| \leq 1.
\end{equation}
With this bound in place, we shift focus to the first term on the right-hand side of \eqref{eq:dummy0}. Note that the matrices $\D_k$ and $\hat{\D}_k$ are both symmetric and positive definite and therefore we can write $\|\D_k^{-1/2}\hat{\D}_k^{1/2} \|^2 = \|\D_k^{-1/2}\hat{\D}_k \D_k^{-1/2}\|$. Notice that the matrix $\D_k^{-1/2}\hat{\D}_k \D_k^{-1/2}$ is block diagonal where its $i$-th diagonal block is given by 
\begin{equation}
{\bf I}+\nabla_{\!\y^i\y^i}g(\y^i_k, \y^j_k;t_k)^{-\frac{1}{2}}\nabla_{\!\y^i\y^i}f^i(\y^i_k;t_k)\nabla_{\!\y^i\y^i}g(\y^i_k, \y^j_k;t_k)^{-\frac{1}{2}}. 
\end{equation}
Using the bounds in Assumptions~\ref{as:first}-\ref{as:2}, and the fact that for positive definite matrices $\lambda_{\min}(\A\mathbold{B}) \geq \lambda_{\min}(\A)\lambda_{\min}(\mathbold{B})$, we obtain that the eigenvalues of the matrix $\D_k^{-1/2}\hat{\D}_k \D_k^{-1/2}$ blocks are bounded below by $1+{2m}/{L}$. Thus, we obtain 
\begin{equation}\label{eq:dummy3}
\|\hat{\D}_k^{-1}\D_k \| =  \|\D_k^{-1/2}\hat{\D}_k \D_k^{-1/2}\| \geq 1+({2m}/{L} ) \; .
\end{equation}
Since the eigenvalues of $\D_k^{-1/2}\hat{\D}_k \D_k^{-1/2}$ are lower bounded by $1+{2m}/L$ we obtain that
%
%\textcolor{blue}{This seems false. I don't think it's possible that the inequality would reverse directions here...}
\begin{equation}\label{eq:dummy31}
\|\D_k^{-1}\hat{\D}_k \| =  \|\D_k^{1/2}\hat{\D}_k^{-1} \D_k^{1/2}\|\leq \left(1+({2m}/{L})\right)^{-1}.
\end{equation}
Substituting the upper bounds in \eqref{eq:dummy21} and \eqref{eq:dummy31} into~\eqref{eq:dummy0} yields 
\begin{equation}\label{eq:varrho}
\|\D_k^{-1/2} \B_k \D_k^{-1/2}\| \leq \varrho \;.
\end{equation}
Note that \eqref{eq:varrho} implies that the eigenvalues of $\D_k^{-1/2} \B_k \D_k^{-1/2}$ are strictly less than one and thus the expansion in \eqref{Hessian_inverse} is valid. 

We use the result in \eqref{eq:varrho} to prove the claim in \eqref{def.H}. Given the approximation in~\eqref{eq:approx}, we know that
\begin{align}\label{eq:finite_geometric}
\|\appH\| 
&\leq \|\D_{k}^{-1}\| \sum_{\tau=0}^K \|\D_k^{-1/2} \B_k \D_k^{-1/2}\|^\tau \nonumber \\
&\leq \frac{1}{m\!+\!\ell/2} \sum_{\tau = 0}^K \varrho^\tau = \frac{1-\varrho^{K+1}}{(m + \ell/2)(1-\varrho)},
\end{align}
where the second inequality comes from formula for a finite geometric series. Moreover, we can derive an upper bound for the RHS of \eqref{eq:finite_geometric}, and use  the definition of $\varrho$ to obtain
\begin{align}\label{eq:finite_geometric1}
\frac{1-\varrho^{K\!+\!1}}{(m\! +\! \ell/2)(1\!-\!\varrho)}\! \leq \! \frac{1}{(m \!+\! \ell/2)(1\!-\!\varrho)}\! =\!  \frac{2m \!+\! L}{m(2m\! +\! \ell)}\! =: \!H
\end{align}
Combining the inequalities in \eqref{eq:finite_geometric} and \eqref{eq:finite_geometric1} the claim in \eqref{def.H} follows. Moreover, the bound on the error $e_k$ follows from~\cite[Proposition~3]{mokhtari2015network1} with the definition of $\varrho$ [cf. \eqref{eq:varrho}], yielding
\begin{equation}
e_k = \|{\bf I} \!-\! \nabla_{\y\y}F(\y_k; t_k) {\H}_{k, (K)}^{-1}\| = \|\D_k^{-\frac{1}{2}} \B_k \D_k^{-\frac{1}{2}}\|^{K+1}\!,
\end{equation}
from which the bound~\eqref{eq.errorE} follows. \qed
%
%\textcolor{red}{[Say more words here, maybe?]}
%
%

\section{Proof of Theorem \ref{theorem.gradient}: case DPC-G}\label{app.DPC-Gconvg}

%\subsection{Case $\booa = 0$, i.e., DPC-G}

%To prove Theorem \ref{theorem.gradient}, we bound the error in the prediction step by the terms that depend on the functional smoothness and the discretization error using Taylor expansions. Then we bound the tracking error of the gradient step using convergence properties of the gradient on strongly convex functions. By substituting the error of the correction step into the prediction step, we establish the main result.

%First, we establish that the criteria in which the bound stated in \eqref{result1} may be satisfied, i.e. that $\rho \sigma$ can be made strictly less than 1, for a sufficiently small $h$. Note that, since $m>0$, there is always a $\gamma$ that makes $\rho < 1$ since $\gamma < m/L^2$. In particular, if we take $\gamma = m/ 2 L^2$, then $\rho$ is minimized and it is equal to 
%
%\begin{equation}\label{rho}
%\rho^* = \sqrt{ 1 - \frac{m^2}{4 L^2}}. 
%\end{equation}
%
%The coefficient $\sigma$ is a monotonically increasing function in $h$, and it is equal to $1$ when $h = 0$; by continuity of $\sigma$ in $h$, we can always find a small enough step-size $h$, such that $\sigma < 1/\rho^*$, which means that we can find a step-size $\gamma$ that yields $\rho \sigma <1$.

First, we establish that discrete-time sampling error bound stated in \eqref{result2} is achieved by the updates of DPC-G. For simplicity, we modify the notation to omit the arguments $\y_k$ and $t_k$ of the function $F$. In particular, define 
\begin{align}\label{eq:thm1_defs}
&\nabla_{\y\y}F := \nabla_{\y\y}F(\y_k; t_k)\; ,\quad \nabla_{t\y} F := \nabla_{t\y} F(\y_k; t_k)\;, \\ 
&\nabla_{\y\y}F^* := \nabla_{\y\y}F({\y}^*(t_{k}); t_k)\;, \ \nabla_{t\y} F^* := \nabla_{t\y} F({\y}^*(t_{k}); t_{k}) . \nonumber
\end{align}
Begin by considering the update of DPC-G, the prediction step, evaluated at a generic point $\y_k$ sampled at the current sample time $t_k$ and with associated optimizer $\y^*(t)$, 
\begin{equation}\label{eq:null_residual}
{\y}^*(t_{k+1}) = {\y}^*(t_{k}) - h\,[\nabla_{\y\y}F^*]^{-1} \nabla_{t\y} F^* + \mathbold{\Delta}_k.
\end{equation}
Rewrite the approximate prediction step $\y_{k+1|k} = \y_k + h\, \p_{k} $ by adding and subtracting the exact prediction step $h\,[\nabla_{\y\y}F]^{-1} \nabla_{t\y} F$, yielding
\begin{multline}\label{eq.ykk1}
{\y}_{k+1|k} = \y_{k} +h \p_{k, (K)} \\ + h\,[\nabla_{\y\y}F]^{-1} \nabla_{t\y} F - h\,[\nabla_{\y\y}F]^{-1} \nabla_{t\y} F.  
\end{multline}
Subtract \eqref{eq:null_residual} from \eqref{eq.ykk1}, take the norm, and apply the triangle inequality to the resulted expression to obtain 
\begin{align}\label{eq.imp}
\|&{\y}_{k+1|k} -  \y^*(t_{k+1})\|   \\  &\leq \|{\y}_{k}\! -\! \y^*(t_{k})\|  \!+
\! h\,\left\|[\nabla_{\y\y}F]^{-1}\nabla_{t\y} F \!-\! [\nabla_{\y\y}F^*]^{-1}\nabla_{t\y} F^* \! \right\| \nonumber \\ 
& \quad + h\,\left\|\p_{k,(K)} - [\nabla_{\y\y}F]^{-1}\nabla_{t\y} F \right\| + \|\mathbold{\Delta}_k\| .\nonumber 
\end{align}
We proceed to analyze the three terms on the RHS of \eqref{eq.imp}. The the last term $ \|\mathbold{\Delta}_k\|$ is bounded above by $h^2 \Delta$ as in \eqref{prop_claim_err_bound}. 

%Applying this substitution into \eqref{eq.imp} implies 
%
%\begin{align}\label{eq.imp00}
%\|{\y}_{k+1|k}& -  \y^*(t_{k+1})\| \leq \|{\y}_{k} - \y^*(t_{k})\|  \\  & + h\,\left\|[\nabla_{\y\y}F]^{-1}\nabla_{t\y} F - [\nabla_{\y\y}F^*]^{-1}\nabla_{t\y} F^*  \right\| \nonumber \\ & + h\,\left\|\d_{k,(K)} - [\nabla_{\y\y}F]^{-1}\nabla_{t\y} F \right\| + h^2 \Delta .\nonumber 
%\end{align}
%
%Set the second term in the right-hand side of \eqref{eq.imp} and the third term as  
%\begin{subequations} \label{eq.ab}
%\begin{align}
%a & = h\,\left\|[\nabla_{\y\y}F]^{-1}\nabla_{t\y} F - [\nabla_{\y\y}F^*]^{-1}\nabla_{t\y} F^*  \right\|, \\
%b & = h\,\left\|\d_{k,(K)} - [\nabla_{\y\y}F]^{-1}\nabla_{t\y} F \right\|.
%\end{align}
%\end{subequations}
%
%We proceed to find an upper bound for the norm terms $a$ and $b$.  
We proceed to find an upper bound for the second summand in the RHS of \eqref{eq.imp}. We use the same reasoning as in~\cite[Appendix~B]{Paper1}, which yields [cf. Eq. (62) of \cite{Paper1}]
\begin{align}\label{eq.dummy20}
h\,\big\|[&\nabla_{\y\y}F]^{-1}\nabla_{t\y} F
  - [\nabla_{\y\y}F^*]^{-1}\nabla_{t\y} F^*  \big\| \\
    &\qquad\quad\leq \frac{C_0 C_1 h}{m^2}\|\y_k - \y^*(t_k)\| + \frac{C_2 h}{m}\|\y_k - \y^*(t_k)\|. \nonumber
\end{align}

Finally, we proceed to analyze the third term in~\eqref{eq.imp}. Rewrite this term using the definition of the prediction step $\p_{k,(K)} =- \appH \nabla_{t\y} F$, and apply the mixed first-order partial derivative bound $\|\nabla_{t\y} F(\y; t)\|\!\leq \!C_0$ stated in Assumption~\ref{as:2} to obtain
\begin{equation}\label{eq.dummy21}
h\,\| \appH\! \nabla_{t\y} F - \nabla_{\y\y}F^{-1} \nabla_{t\y} F \| \!\leq\! C_0 h\,\|\appH \!-  \nabla_{\y\y}F^{-1}\| \; .
\end{equation}
Observe that $\|\appH \!-  \nabla_{\y\y}F^{-1}\|$ is bounded above by $\|\nabla_{\y\y}F^{-1}\|\|\nabla_{\y\y}F \appH - {\bf I}\|$. This observation in conjunction with the upper bound for the error vector $e_k=\|\nabla_{\y\y}F \appH - {\bf I}\|$ in Proposition \ref{th.hessian} implies that
\begin{align}\label{eq.dummy22}
&\|\appH -  \nabla_{\y\y}F^{-1}\| \leq   \\
 &  \qquad \leq \|\nabla_{\y\y}F^{-1}\|\|\nabla_{\y\y}F \appH - {\bf I}\|  
 \leq \frac{\varrho^{K+1}}{m} .\nonumber
\end{align}
Combining the results in \eqref{eq.dummy21} and \eqref{eq.dummy22} shows that the third in the RHS of \eqref{eq.imp} is upper bounded by
\begin{align}\label{eq.dummy33}
h\,\left\|\p_{k,(K)} - [\nabla_{\y\y}F]^{-1}\nabla_{t\y} F \right\| 
 \leq \frac{hC_0}{m} \varrho^{K+1}\ .
\end{align}
By substituting the bounds in~\eqref{eq.dummy20} and \eqref{eq.dummy33} into~\eqref{eq.imp} and considering the definitions of $\sigma$ in~\eqref{def_sigma} and $\Gamma(\varrho, K)$ in Theorem \ref{theorem.gradient} we obtain 
\begin{equation}\label{eq.result10}
\|{\y}_{k+1|k} \!-\! \y^*(t_{k+1})\|\! \leq \! \sigma \|{\y}_{k} \!-\! \y^*(t_{k})\| \!+\! h\,\Gamma(\varrho, K)\! + \!h^2 \Delta.
\end{equation}

For the correction step [cf. \eqref{correctionstep}] , we may use the standard property of projected gradient descent for strongly convex functions with Lipschitz gradients. The Euclidean error norm of the projected gradient descent method converges linearly as
\begin{equation}\label{eq.fin2}
\|\y_{k+1} - \y^*(t_{k+1})\| \leq \rho \|\y_{k+1|k } - \y^*(t_{k+1})\| .
\end{equation}
where $\rho = \max\{|1-\gamma m|,|1- \gamma (L+M)|\}$; see e.g., \cite{Paper1} or~\cite{Ryu2015}. Plug the correction error in \eqref{eq.fin2} into the prediction error in \eqref{eq.result10} to obtain
\begin{equation}\label{eq.fin_t_k}
\|{\y}_{k+1} - \y^*(t_{k+1})\| \leq \rho \sigma \|{\y}_{k} - \y^*(t_{k})\| +  \rho\, \varphi,
\end{equation}
where $\varphi:= h\, \Gamma(\varrho, K) + h^2 \Delta$. Therefore,%Now recursively apply the relationship \eqref{eq.fin_t_k} backwards in time to the initial time sample to write 
\begin{equation}\label{eq.fin_t_00}
\|{\y}_{k+1} - \y^*(t_{k+1})\| \leq (\rho \sigma)^{k+1} \|{\y}_{0} - \y^*(t_{0})\| + \rho \,\varphi \sum_{i=0}^k (\rho \sigma)^i.
\end{equation}
Substitute $k+1$ by $k$ and simplify the sum in \eqref{eq.fin_t_00}, making use of the fact that $\rho \sigma < 1$, which yields
\begin{equation}\label{eq.fin_t_0000}
\|{\y}_{k}\! - \!\y^*(t_{k})\| \!\leq\! (\rho \sigma)^{k} \|{\y}_{0} - \y^*(t_{0})\| \!+ \!\rho \,\varphi \!\left[\! \frac{1-(\rho\sigma)^k}{1-\rho\sigma}\!\right].
\end{equation}
Observing \eqref{eq.fin_t_0000} together with the definition of $\varphi$, \eqref{result2} follows. In particular, if $\rho \sigma < 1$ (that is~\eqref{eq.hh} holds), the sequence $\{\y_{k}\}$ converges Q-linearly to $\y^*$ up to an error bound as  
\begin{align}\label{eq.fin_t_0000bis}
&\limsup_{k\to \infty}  \|{\y}_{k} - \y^*(t_{k})\| = \\ &\quad \left[\! \frac{\rho}{1-\rho\sigma}\!\right] \!\left(h\, \Gamma(\varrho, K) + h^2 \left[\frac{C_0 C_2}{m^2} + \frac{C_3}{2m}+ \frac{C_0^2C_1}{2m^3} \right]\right) .\nonumber
\end{align}

To establish the result stated in \eqref{result1}, observe that in the worst case, we may upper bound the term $\|[\nabla_{\y\y}F]^{-1}\nabla_{t\y} F - [\nabla_{\y\y}F^*]^{-1}\nabla_{t\y} F^* \|$ in \eqref{eq.imp} by $\left\|[\nabla_{\y\y}F]^{-1}\nabla_{t\y}F\|+\| [\nabla_{\y\y}F^*]^{-1}\nabla_{t\y} F^*  \right\|  $ which yields
\begin{align}\label{eq.worst_case}
\left\|[\nabla_{\y\y}F]^{-1}\nabla_{t\y}F - [\nabla_{\y\y}F^*]^{-1}\nabla_{t\y} F^*  \right\| 
 & \leq \frac{2 C_0}{m}.
\end{align}
Substituting the upper bound in~\eqref{eq.worst_case} into \eqref{eq.imp} yields
\begin{align}\label{eq.a_subst}
\|{\y}_{k+1|k} - \y^*(t_{k+1})\| \leq \|{\y}_{k} - \y^*(t_{k})\| + h\, \frac{2 C_0 }{m} + \varphi.
\end{align}
%
%To simplify the notation we define a new constant $2 h { C_0 }/{m} := 2 h { C_0 }/{m} $ and 
%
Using the definition $\varphi=h\, \Gamma(\varrho, K) + h^2 \Delta$ and observing the relation in \eqref{eq.fin2}, we can write 
\begin{equation}\label{eq.sampling_error2}
\|{\y}_{k+1} - \y^*(t_{k+1})\| \leq \rho \|{\y}_{k} - \y^*(t_{k})\| + \rho \left[2 h \frac{ C_0 }{m}  + \varphi\right].
\end{equation}
Recursively applying \eqref{eq.sampling_error2} backwards to the initial time% and use the analogous logic from \eqref{eq.fin_t_k} to \eqref{eq.fin_t_0000} to write 
\begin{equation}\label{eq.rho_error_dependence}
\!\!\|{\y}_{k} \!-\! \y^*(t_{k})\| \!\leq \!\rho^{k} \|{\y}_{0} \!- \!\y^*(t_{0})\|\! + \! \rho \left[\! \frac{ 2 hC_0 }{m} \! +\! \varphi\!\right] \!\left[ \frac{1-\rho^{k}}{1-\rho}\right],
\end{equation}
and by sending $k\to \infty$, we can simplify \eqref{eq.rho_error_dependence} as
\begin{align}\label{eq.fin_t_0000tris}
\limsup_{k\to \infty} & \|{\y}_{k} - \y^*(t_{k})\|=  \frac{2 C_0\rho h}{m(1-\rho)} +  \\ &  \frac{\rho}{1-\rho}  h \Gamma(\varrho, K) + h^2 \frac{\rho}{1-\rho} \left[\frac{C_0 C_2}{m^2} + \frac{C_3}{2m}+ \frac{C_0^2C_1}{2m^3} \right],\nonumber
\end{align}
which is \eqref{result1}. The result in \eqref{eq.fin_t_0000tris} holds if $\rho < 1$, which is the case if $\gamma < 2/(L+M)$, as stated in Theorem \ref{theorem.gradient}. %Considering the definition of $\rho$ in \eqref{def_sigma} we require 
%
%\begin{equation}\label{sososo}
%\rho := \max\{|1-\gamma m|,|1-\gamma (L+M)|\} < 1.
%\end{equation}
%
%Solving \eqref{sososo} for $\gamma$ and recalling that $m \leq L+M$ by Assumptions~\ref{as:first} and \ref{as:last}, we obtain that the step-size $\gamma$ should satisfy $\gamma < 2/(L+M)$, as stated in Theorem \ref{theorem.gradient}.
%
\qed
%
% APPROXIMATE TIME !
%
%
%\subsection{Case $\booa = 1$, i.e. DAPC-G}
%\begin{comment}
%%PUT THIS ON ARXIV
\section{Proof of Theorem \ref{theorem.gradient}: case DAPC-G} \label{app.DAPC-Gconvg}

We establish the tracking performance of the DAPC-G method by following a similar line of reasoning as that which yields the DPC-G error bounds. To do so, we characterize the error coming from the approximate time derivative in~\eqref{fobd}. In particular, consider the Taylor expansion of the gradient ${\nabla}_{\y}{F}(\y_{k}; t_{k-1})$ near the point $(\y_k,t_k)$ which is given by
%%%
\begin{align}\label{taylor_series100}
{\nabla}_{\y}{F}(\y_{k}; t_{k-1}) 
&= {\nabla}_{\y}{F}(\y_k; t_k) - h\,{\nabla}_{t\y}{F}(\y_k; t_k) \nonumber \\
&\quad+  \frac{h^2}{2}  {\nabla}_{tt\y}{F}(\y_k; s) .
\end{align}
for a particular $s\!\in\![t_{k-1}, t_k]$. Regroup terms in \eqref{taylor_series100} to obtain 
%%%
\begin{align}\label{taylor_series200}
{\nabla}_{t\y}{F}(\y_k;  t_k) & = \frac{{\nabla}_{\y}{F}(\y_k; t_k)\! -\! {\nabla}_{\y}{F}(\y_{k}; t_{k-1})}{h}\!\nonumber \\
&\quad+ \!  \frac{h}{2}  {\nabla}_{tt\y}{F}(\y_k; s) \; .
\end{align}
%%%%
Use the definition of the approximate partial mixed gradient $\tilde{\nabla}_{t\y}F(\y_k; t_k) $ in~\eqref{fobd} and the expression for the exact mixed gradient $ {\nabla}_{\y}{F}(\y_k; t_k) $ in \eqref{taylor_series200} to obtain 
%%%
\begin{equation}\label{taylor_series300}
{\nabla}_{t\y}{F}(\y_{k}; t_{k}) -\tilde{\nabla}_{t\y}F(\y_k; t_k) =\frac{h}{2}  {\nabla}_{tt\y}{F}(\y_k; s).
\end{equation}
%%%
Based on Assumption \ref{as:last} the norm ${\nabla}_{tt\y}{F}(\y_k; s)$ is bounded above by a constant $C_3$, which yields the upper estimate on the error of the partial mixed gradient approximation [cf. also Eq. (100) of \cite{Paper1}]
%%%
\begin{equation}\label{taylor_series400}
\|{\nabla}_{t\y}{F}(\y_{k}; t_{k}) -\tilde{\nabla}_{t\y}F(\y_k; t_k)\| \leq \frac{hC_3}{2}.
\end{equation}
Consider the approximate prediction step of the DAPC-G algorithm. 
By adding and subtracting the prediction direction $\p_{k,(K)}$ to the right-hand side of the update in \eqref{eq.ykk1} we obtain
%%% 
\begin{align}\label{new_number}
\y_{k+1|k} &= \y_k + h \,\p_{k,(K)} \\ 
&\quad 
+h\, \appH \!\! \left({\nabla}_{t\y}{F}(\y_k; t_k)  -\tilde{\nabla}_{t\y}F(\y_k; t_k)\right)\! . \nonumber 
\end{align}
Proceed exactly as the proof of DPC-G in Theorem \ref{theorem.gradient}, starting from \eqref{eq:null_residual}, and observe the presence of the extra term due to the time-derivative approximation error
\begin{align}\label{new_number}
\hskip-0.2cm h \| \appH  \!\! \left({\nabla}_{t\y}{F}(\y_k; t_k)  -\tilde{\nabla}_{t\y}F(\y_k; t_k)\right)\| \leq \frac{h^2}{2}  {C_3H}
\end{align}
in the prediction errors stated in~\eqref{eq.result10} and \eqref{eq.a_subst}. This approximation error term is obtained by combing the error estimate in~\eqref{taylor_series400} with the fact that the eigenvalues of the approximated Hessian $\appH$ are upper bounded as stated in~\eqref{def.H} (Proposition~\ref{th.hessian}). Incorporating this extra error term from time-derivative approximation, the same logic from the proof of Theorem \ref{theorem.gradient} (case DPC-G) yields that if $\rho \sigma < 1$, the sequence $\{\y_{k}\}$ converges Q-linearly to $\y^*$ up to an error bound as  
\begin{align}\label{eq.fin_t_00000bis}
&\limsup_{k\to \infty}  \|{\y}_{k} - \y^*(t_{k})\| = 
\\ &  \left[\! \frac{\rho}{1\!-\!\rho\sigma}\!\right] \!\left[h\, \Gamma(\varrho, K) \!+\! h^2 \left[\frac{C_0 C_2}{m^2}\! +\! \frac{C_3}{2m}(1\!+\!m H)\!+\! \frac{C_0^2C_1}{2m^3} \right]\right].\nonumber
\end{align}
If $\rho < 1$, the sequence $\{\y_{k}\}$ converges Q-linearly to $\y^*$ up to an error bound as
\begin{align}\label{eq.fin_t_00000tris}
\limsup_{k\to \infty} & \|{\y}_{k} - \y^*(t_{k})\|\! = h\!\left[\! \frac{2\rho C_0}{m(1-\rho)}\!\right]\! + \frac{\rho}{1-\rho} h \Gamma(\varrho, K) \nonumber \\ & \hskip-1.5cm + h^2 \left[\frac{C_0 C_2}{m^2} + \frac{C_3}{2m}(1+m H)+ \frac{C_0^2C_1}{2m^3} \right]\!\left[\! \frac{\rho}{1-\rho}\!\right].
\end{align}
%%
%\end{comment}

%%%%%%%%%%%%%%%%%%%%%%%%%%%%%%%%%
%%%%%%%%%%%%%%%%%%%%%%%%%%%%%%%%%
%%%%%%   S  E  C  T  I  O  N    %%%%%%%%%%%%%%%%
%%%%%%%%%%%%%%%%%%%%%%%%%%%%%%%%%
%%%%%%%%%%%%%%%%%%%%%%%%%%%%%%%%%
\section{Proof of Theorem~\ref{th.newton}: case DPC-N}\label{newton.proof1}

%\subsection{Case $\booa = 0$, i.e., DPC-N}

Since DPC-N and DPC-G are identical in their prediction steps, we may consider the prediction error result established during the proof Theorem~\ref{theorem.gradient}, i.e. the expression in \eqref{eq.result10} with $k = 0$. We turn our attention to the correction step, and consider in particular the gap to the optimal trajectory before and after correction at time $t_{k+1}$ as
\begin{equation}\label{eq:optimality_gap}
\|\y_{k+1} - \y^*(t_{k+1})\| =  \|\y_{k+1|k} -\gamma\,{\bf H}_{k+1|k}^{-1} \nabla_{\y}F- \y^*(t_{k+1})\|.
\end{equation}
Subsequently, we simplify notation by defining the shorthands
\begin{align}\label{eq:thm2_defs}
&\nabla_{\y}F \! : = \!\nabla_{\y}F(\y_{k+1|k}; t_{k+1}) \; ,
\nabla_{\y\y}F \!:=\! \nabla_{\y\y}F(\y_{k+1|k}; t_{k+1})\;, \nonumber \\ 
&\nabla_{\y}F^* \!\!:=\! \nabla_{\y}F({\y}^*\!(t_{k+1}); \!t_{k+1}\!) ,\!
 \nabla_{\y\y}F^*\!\! \!:=\!\! \nabla_{\y\y}\!F({\y}^*\!(t_{k+1}\!); \!t_{k+1}\!),  \nonumber\\
& \nabla_{t\y} F := \!\!\nabla_{t\y} \!F(\y_{k+1|k}; \!t_{k+1}\!) \; ,
 \nabla_{t\y} F^*\!\! :=\!\! \nabla_{t\y} \!F({\y}^*(\!t_{k+1}\!);\! t_{k+1}) \; . 
\end{align}
Add and subtract $\gamma\,\nabla_{\y\y}F^{-1}\nabla_{\y} F$ to the expression inside the norm in \eqref{eq:optimality_gap} which is the exact damped Newton, and use the the triangle inequality to obtain
\begin{align}\label{eq:optimality_gaptris}
\|\y_{k+1} \!-\! \y^*(t_{k+1})\| &\!\leq \!
 \|\y_{k+1|k} -\!\!\gamma \nabla_{\y\y}F^{-1}\nabla_{\y} F \!- \!\y^*(t_{k+1})\|\nonumber \\ 
& \ + \gamma\|(\nabla_{\y\y}F^{-1}-{\bf H}_{k+1|k}^{-1})\nabla_{\y} F\|. 
\end{align}
We proceed to bound the two terms in the RHS of~\eqref{eq:optimality_gaptris}. Consider the first term: left multiply by $\nabla_{\y\y}F$ and its inverse, and left factor out the Hessian inverse $\nabla_{\y\y}F^{-1}$. Making use of the Cauchy- Schwartz inequality, the first term of right-hand side of~\eqref{eq:optimality_gaptris} is bounded above as
%%%
\begin{multline}\label{eq.aryan3}
 \|\y_{k+1|k} -  \gamma \nabla_{\y\y}F^{-1}\nabla_{\y}F- \y^*(t_{k+1})\| \leq \\ (1-\gamma)\|\y_{k+1|k}- \y^*(t_{k+1})\|+ \\
\gamma\|\nabla_{\y\y}F^{-1}\|\|\nabla_{\y\y}F(\y_{k+1|k} \!-\! \y^*(t_{k+1})) \!-\! \nabla_{\y}F\| \; .%\nonumber
\end{multline}
%%%
We use now the same arguments as in~\cite[Appendix~C, Eq.s (83)-(85)]{Paper1} to show that~\eqref{eq.aryan3} can be upper bounded by
\begin{align}\label{eq.aryan3000}
%\left\|\y_{k+1|k} -  \nabla_{\y\y}F^{-1}\nabla_{\y}F- \y^*(t_{k+1})\right\|
(1-\gamma)\|  \y_{k+1|k} - \y^*(t_{k+1})\|+  \frac{\gamma C_1}{2m} \|  \y_{k+1|k} - \y^*(t_{k+1})\|^2 \; .
\end{align}
%%%
With this upper estimate in place for the first term on the right hand side of~\eqref{eq:optimality_gaptris}, we shift focus to the second term. Use the Cauchy-Schwartz inequality to obtain $ \|(\nabla_{\y\y}F^{-1}-{\bf H}_{k+1|k}^{-1})\nabla_{\y} F\|$ is bounded above by $  \|\nabla_{\y\y}F^{-1}\|\|\nabla_{\y\y}F {\bf H}_{k+1|k}^{-1} - {\bf I}\|  \|\nabla_{\y} F\|$. Use the upper bound $1/m$ for the spectrum of $\nabla_{\y\y}F^{-1}$ and the Hessian approximation error [cf. \eqref{eq.errorE}] to obtain
\begin{align}\label{eq.aryan2t_1}
\|(\nabla_{\y\y}F^{-1}-{\bf H}_{k+1|k}^{-1})\nabla_{\y} F\|\leq \frac{\varrho^{K'+1} }{m}\|\nabla_{\y} F\| .
\end{align}
Now, focusing on the second term in the product on the right-hand side of \eqref{eq.aryan2t_1}, we use of the optimality criterion of $\y^*(t_{k+1})$, which is equivalent to $\nabla_{\y}F^* = \mathbf{0}$, to write
\begin{align}\label{eq.aryan2t_3}
\!\!\!\!\!\!\|\nabla_{\y} F\| \!= \! \|\nabla_{\y} F\!-\!\nabla_{\y}F^*\|\!
\leq\! (L\!+\!M) \|\y_{k+1|k} \!-\! \y^*(t_{k+1}) \|,
\end{align}
where we have used the Lipschitz property of the gradients. Substituting the upper bound in \eqref{eq.aryan2t_3} into \eqref{eq.aryan2t_1} and considering the definition $\Gamma(\varrho, K')=(C_0/m){\varrho^{K'+1} }$ lead to   
\begin{align}\label{eq.aryan2t_4}
&\|(\nabla_{\y\y}F^{-1}-{\bf H}_{k+1|k}^{-1})\nabla_{\y} F\| 
\nonumber\\
&\qquad\qquad\leq\frac{L+M}{C_0} \Gamma(\varrho, K') \|\y_{k+1|k} - \y(t_{k+1})\| .
\end{align}
Apply the bounds \eqref{eq.aryan3000} - \eqref{eq.aryan2t_4} to the right-hand side of \eqref{eq:optimality_gaptris}
\begin{align}\label{eq.dummy300}
\|\y_{k+1} &- \y^*(t_{k+1})\| \leq  \gamma \frac{C_1}{2 m} \|\y_{k+1|k} - \y^*(t_{k+1})\|^2  \\
&\quad + \Big(\gamma\frac{L+M}{C_0} \Gamma(\varrho, K') + 1 - \gamma\Big)\|\y_{k+1|k} - \y^*(t_{k+1})\| \; . \nonumber
\end{align}
%%%
Now we consider the prediction step, which by~\eqref{eq.result10} we have
\begin{equation}\label{eq.result10_another}
\|{\y}_{k+1|k} - \y^*(t_{k+1})\| \leq \sigma \|{\y}_{k} - \y^*(t_{k})\| + \varphi,
\end{equation}
with $\varphi=h\, \Gamma(\varrho, K) + h^2 \Delta$ as defined in Appendix \ref{app.DPC-Gconvg}. Substituting the relation \eqref{eq.result10_another} into \eqref{eq.dummy300} allows us to write
\begin{align}\label{eq.dummy300_another}
\|\y_{k+1} - \y^*&(t_{k+1})\| \leq  \gamma \frac{C_1}{2 m} 
 \left(\sigma \|\y_{k} - \y^*(t_{k})\|+ \varphi\right)^2 \\
&\!\!\!\!\!\!\!\!\!\!\!\!\!\!+\!  \Big(\gamma\frac{L+M}{C_0} \Gamma(\varrho, K') + 1 - \gamma\Big)
 \left(\sigma \|{\y}_{k}\! -\! \y^*(t_{k})\| \!+\! \varphi\right).\nonumber
\end{align}

The right-hand side of \eqref{eq.dummy300_another} is a quadratic function of the error $\|\y_{k} - \y^*(t_{k})\|$ at time $t_k$, which upper bounds the error sequence at the subsequent time $t_{k+1}$. For certain selections of parameters $K$, $K'$, and $h$, \eqref{eq.dummy300_another} defines a contraction. To determine the conditions for which this occurs, we solve for an appropriate radius of contraction. Let $\tau>0$ be a positive scalar such that
\begin{align}\label{q_inequality}
\alpha_2 \|\y_{k} - \y^*(t_{k})\|^2 + &\alpha_1 \|\y_{k} - \y^*(t_{k})\| + \alpha_0  \\
&\leq \tau \|\y_{k} - \y^*(t_{k})\| + \alpha_0 \; , \nonumber
\end{align}
%
%Call $q_{k+1} = \|\y_{k+1} - \y^*(t_{k+1})\|$ and $q_k = \|\y_{k} - \y^*(t_{k})\|$. The relation~\eqref{eq.dummy300_another} implies
%\begin{equation}\label{eq.dummy301_another}
%q_{k+1} \leq  \alpha_2 q_k^2 + \alpha_1 q_k + \alpha_0,
%\end{equation}
%with
where the coefficients $\alpha_0$, $\alpha_1$, and $\alpha_2$ of the quadratic polynomial of the error $\|\y_{k} - \y^*(t_{k})\|$ are defined from the right-hand side of \eqref{eq.dummy300_another}, and given as
\begin{align}\label{def_alpha1}
\alpha_2 &= \gamma \frac{C_1}{2 m} \sigma^2,
 \, \alpha_1 = \sigma\left[\gamma\frac{C_1}{m} \varphi + \gamma\frac{L+M}{C_0}\Gamma(\varrho, K') + 1 - \gamma\right], \nonumber \\
\alpha_0 &= \varphi\left[\gamma\frac{C_1}{2 m} \varphi + \gamma\frac{L+M}{C_0}\Gamma(\varrho, K')+1-\gamma\right].
\end{align}
Based on \eqref{q_inequality}, to guarantee the Q-linear convergence of the error sequence $\|\y_{k} - \y^*(t_{k})\|$, we require $\tau<1$, which by simple algebra it is satisfied if 
\begin{equation}\label{eq:coefficient_conditions}
\alpha_1 < \tau, \quad \|\y_{0} - \y^*(t_{0})\| \leq (\tau-\alpha_1)/\alpha_2 \; .
\end{equation}
Finally, we need to require that the second condition in~\eqref{eq:coefficient_conditions} holds true for all $k$, that is $\|\y_{k+1} - \y^*(t_{k+1})\|\leq \|\y_{k} - \y^*(t_{k})\|$, which implies
\begin{equation}\label{eq:coefficient_conditions2}
\tau (\tau-\alpha_1)/\alpha_2 +\alpha_0 \leq (\tau-\alpha_1)/\alpha_2.
\end{equation}
Conditions~\eqref{eq:coefficient_conditions}-\eqref{eq:coefficient_conditions2}, the definitions of $\alpha_0, \alpha_1$, and $\alpha_2$ in~\eqref{def_alpha1}, and $\varphi=h\, \Gamma(\varrho, K) + h^2 \Delta$ establish the small enough conditions on sampling period and optimality gap as well as the large enough conditions on the approximation levels $K, K'$ in Theorem \ref{th.newton}, for any chosen $\tau < 1$.   
In particular, the terms $\alpha_0$ and $\alpha_1$ are polynomial functions of the sampling period $h$ and the approximation levels $K$ and $K'$. Conditions~\eqref{eq:coefficient_conditions}-\eqref{eq:coefficient_conditions2} describe a system of nonlinear inequalities for any fixed $1-\gamma <\tau < 1$. For arbitrarily small $h$ and large $K$ and $K'$, $\alpha_0$ and $\alpha_1$ can be made $0$ and $1-\gamma$, respectively. Formally
\begin{equation}
\lim_{h\to 0, \,K,K' \to \infty} \alpha_0 = 0, \quad  \lim_{h\to 0, \,K,K' \to \infty} \alpha_1 = 1-\gamma.
\end{equation}
When $\alpha_0 = 0$ and $\alpha_1 = 1-\gamma$ conditions~\eqref{eq:coefficient_conditions}-\eqref{eq:coefficient_conditions2} hold with attraction region $\bar{R} =  {2 m} (\tau-1+\gamma)/{\gamma C_1 \sigma^2}$. Since $\alpha_0$ and $\alpha_1$ go their limits monotonically with $h, K, K'$, then -- by continuity -- there exists a large enough $\bar{K}$ and small enough attraction region $\bar{R}$ for which conditions \eqref{eq:coefficient_conditions}-\eqref{eq:coefficient_conditions2} can be achieved. %The condition  $\alpha_1 < \tau$ together with the result in~\eqref{result_37} display the trade-off between high converge rates, i.e., small $\tau$'s, and increased computational and communication complexities (small $\tau$'s implies small $h$ and big $K$ and $K'$). 
The convergence region in~\eqref{eq:coefficient_conditions} is dictated by the Newton step and by the step-size choice $\gamma$. If the cost function is a time-varying quadratic function, then $C_1 = 0$, and the convergence region is the whole space. If $\gamma$ is chosen very small, then also in this case the convergence region becomes arbitrarily large, as expected.  

By selecting the error polynomial coefficients as \eqref{def_alpha1} and recursively applying~\eqref{q_inequality} backward in time, we obtain
\begin{equation}\label{q_inequality1}
\|\y_{k} - \y^*(t_{k})\| \leq \tau^k \|\y_{0} - \y^*(t_{0})\| +  \alpha_0 \left[\frac{1-\tau^k}{1-\tau}\right] \; ,
\end{equation}
which since $\tau < 1$, the right-side of the \eqref{q_inequality1} is finite, implying~\eqref{result_37}, after the expansion of the coefficients. \qed
%

%%%%%%%%%%%%%%%%%%%%%%%%%%%%%%%%%
%%%%%%%%%%%%%%%%%%%%%%%%%%%%%%%%%
%%%%%%   S  E  C  T  I  O  N    %%%%%%%%%%%%%%%%
%%%%%%%%%%%%%%%%%%%%%%%%%%%%%%%%%
%%%%%%%%%%%%%%%%%%%%%%%%%%%%%%%%%

\section{Theorem~\ref{th.newton}: case DAPC-N}\label{newton.proof2}
%\subsection{Case $\booa = 1$, i.e., DAPC-N}
%
%\emph{(Sketch)}
Proceed with analogous logic to that which establishes the convergence of DPC-N in Theorem~\ref{th.newton}. Instead of the error contraction in~\eqref{eq.result10_another} when the exact time-derivative is used, we obtain
\begin{equation}\label{eq.result10_another22}
\|{\y}_{k+1|k} - \y^*(t_{k+1})\| \leq \sigma \|{\y}_{k} - \y^*(t_{k})\| + \varphi + \frac{h^2}{2}\, {C_3\,H}\; ,
\end{equation}
where the last term on the right-hand side comes from the error due to approximating the time derivative, and is derived in \eqref{taylor_series400}. Proceeding in a similar manner as to that which yields a quadratic polynomial of the error sequence $\|{\y}_{k} - \y^*(t_{k})\|$ in the proof of Theorem~\ref{th.newton}, replacing $\varphi$ with $\varphi + h^2\, {C_3} H/2$, and solving for the required conditions on the problem parameters to obtain a contraction, the proof is completed. In particular, consider $\varphi' = \varphi + \frac{h^2}{2}\, {C_3\,H}$ and define the coefficients $\alpha_0'$, $\alpha_1'$, and $\alpha_2$ as
\begin{align}\label{def_alpha1}
\alpha_2 &= \gamma\frac{C_1}{2 m} \sigma^2,
 \,\alpha_1' = \sigma\left[\gamma\frac{C_1}{m} \varphi' + \gamma\frac{L+M}{C_0}\Gamma(\varrho, K') + 1-\gamma\right], \nonumber \\
\alpha_0' &= \varphi'\left[\gamma\frac{C_1}{2 m} \varphi' + \gamma\frac{L+M}{C_0}\Gamma(\varrho, K') + 1-\gamma\right].
\end{align}
If the modified conditions~\eqref{eq:coefficient_conditions}-\eqref{eq:coefficient_conditions2} in terms of $\alpha_0'$, $\alpha_1'$, and $\alpha_2$ hold, then convergence of $\{\y_k\}$ goes as 
\begin{equation}\label{q_inequality1}
\|\y_{k} - \y^*(t_{k})\| \leq \tau^k \|\y_{0} - \y^*(t_{0})\| +  \alpha_0' \left[\frac{1-\tau^k}{1-\tau}\right] \,.
\end{equation}
%
%\end{comment}

%\section{Proof of Theorem~\ref{th.ANT_convg} }\label{ap.ant}

%%%%%%%%%%%%%%%%%%%%%%%%%%%%%%%%%%%
%%%%%%%%%%%%%%%%%%%%%%%%%%%%%%%%%%%
%%%%%     P  R   O  O   F    %%%%%%%%%%%%%%%%%%%%
%%%%%%%%%%%%%%%%%%%%%%%%%%%%%%%%%%%
%%%%%%%%%%%%%%%%%%%%%%%%%%%%%%%%%%%

\bibliographystyle{ieeetr}
\bibliography{PaperCollection2}

\end{document}